\documentclass[12pt,a4paper]{book}
\usepackage{silence}
\usepackage[inner=3.75cm,outer=3.75cm,centering]{geometry}
\usepackage{times}
\usepackage[amsthm]{ntheorem}
\usepackage{amsfonts,amsgen,amssymb,amsmath,amsbsy}
\usepackage{latexsym}
\usepackage[croatian]{babel}
\usepackage[T1]{fontenc}
\usepackage[utf8]{inputenc}
\usepackage[nottoc]{tocbibind}
\usepackage{calc}
\usepackage{ifthen}
\usepackage{mathtools}
\usepackage{listings}
\usepackage{lmodern}
\usepackage{textcomp}
\usepackage{minitoc}
\usepackage[hidelinks]{hyperref}
\usepackage{thmtools}
\hypersetup{
  colorlinks   = true, 
  urlcolor     = blue, 
  linkcolor    = blue, 
  citecolor   = red 
}
\usepackage{etex}
\usepackage{graphics}
\usepackage{tabularx}
\usepackage{pstricks,pst-node,pst-text,pst-3d}
\usepackage{pst-slpe}
\usepackage{pst-plot}
\usepackage{pst-grad}
\usepackage{pst-coil}
\usepackage{pst-3dplot}
\usepackage{pst-circ}
\usepackage{subfigure}
\usepackage{graphicx}
\usepackage[font=small,labelfont=bf,tableposition=top]{caption}
\usepackage{pstricks-add}
\usepackage{tikz}

\usetikzlibrary{calc,trees,positioning,arrows,chains,shapes.geometric,%
decorations.pathreplacing,decorations.pathmorphing,shapes,%
matrix,shapes.symbols,plotmarks,decorations.markings,shadows}

\usepackage{titleps}
\theoremstyle{definition}
\newtheorem{Primjer}{\textit{Primjer}}[chapter]
\newtheorem{Vjezba}{Vje\v zba}[chapter]
\newtheorem{Remark}{Primjedba}[chapter]

\newenvironment{dokaz}
{$\mathbf{Dokaz:}$} {~ \hfill {\scalebox{0.4}{\stechak L}} \vspace{-0.1cm}}

\newenvironment{remark}
{\begin{Remark}}{~ \hfill {\scalebox{0.6}{\stechak 5}} \end{Remark}\vspace{-0.1cm}}

\usepackage[framemethod=TikZ,ntheorem]{mdframed}

\newcounter{teorem}[chapter]
\renewcommand{\theteorem}{\thechapter .\arabic{teorem}}
\newenvironment{teorem}[1][]{%
\refstepcounter{teorem}%
\ifstrempty{#1}%
 {\mdfsetup{%
   frametitle={%
    \tikz[baseline=(current bounding box.east),outer sep=0pt]
    \node[anchor=east,rectangle,rounded corners,drop shadow,fill=black!40!blue]
         {\strut \textcolor{white}{Teorem~\theteorem}};}}
 }%
{\mdfsetup{%
  frametitle={%
   \tikz[baseline=(current bounding box.east),outer sep=0pt]
   \node[anchor=east,rectangle,rounded corners,drop shadow,fill=black!40!blue]
        {\strut \textcolor{white}{Teorem~\theteorem:~#1}};}}%
 }%
\mdfsetup{innertopmargin=10pt,linecolor=blue!60,%
       linewidth=2pt,topline=true,roundcorner=5pt,
       frametitleaboveskip=\dimexpr-\ht\strutbox\relax,}
   \begin{mdframed}[]\relax%
}
{\end{mdframed}}
\newenvironment{lema}[1][]{%
\refstepcounter{teorem}%
\ifstrempty{#1}%
 {\mdfsetup{%
   frametitle={%
    \tikz[baseline=(current bounding box.east),outer sep=0pt]
    \node[anchor=east,rectangle,rounded corners,drop shadow,fill=black!60!green]
         {\strut \textcolor{white}{Lema~\theteorem}};}}
 }%
{\mdfsetup{%
  frametitle={%
   \tikz[baseline=(current bounding box.east),outer sep=0pt]
   \node[anchor=east,rectangle,rounded corners,drop shadow,fill=black!60!green]
        {\strut \textcolor{white}{Lema~\theteorem:~#1}};}}
 }%
\mdfsetup{innertopmargin=10pt,linecolor=black!40!green,%
       linewidth=2pt,topline=true,roundcorner=5pt,
       frametitleaboveskip=\dimexpr-\ht\strutbox\relax,}
   \begin{mdframed}[]\relax%
}
{\end{mdframed}}

\declaretheorem{defn}
\renewenvironment{defn}[1][]{%
\refstepcounter{teorem}%
\ifstrempty{#1}%
 {\mdfsetup{%
   frametitle={%
    \tikz[baseline=(current bounding box.east),outer sep=0pt]
    \node[anchor=east,rectangle,rounded corners,drop shadow,fill=black!40!red]
         {\strut \textcolor{white}{Definicija~\theteorem}};}}
 }
{\mdfsetup{%
  frametitle={%
   \tikz[baseline=(current bounding box.east),outer sep=0pt]
   \node[anchor=east,rectangle,rounded corners,drop shadow,fill=black!40!red]
        {\strut \textcolor{white}{Definicija~\theteorem:~#1}};}}%
 }%
\mdfsetup{innertopmargin=10pt,linecolor=black!20!red,%
       linewidth=2pt,roundcorner=5pt,topline=true,
       frametitleaboveskip=\dimexpr-\ht\strutbox\relax,}
   \begin{mdframed}[]\relax%
}
{\end{mdframed}}

\newenvironment{notacija}[1][]{%
\ifstrempty{#1}%
 {\mdfsetup{%
   frametitle={%
    \tikz[baseline=(current bounding box.east),outer sep=0pt]
    \node[anchor=east,rectangle,rounded corners,drop shadow,fill=black!40!gray]
         {\strut \textcolor{white}{Notacija}};}}
 }
{\mdfsetup{%
  frametitle={%
   \tikz[baseline=(current bounding box.east),outer sep=0pt]
   \node[anchor=east,rectangle,rounded corners,drop shadow,fill=black!40!gray]
        {\strut \textcolor{white}{Notacija}};}}%
 }%
\mdfsetup{innertopmargin=10pt,linecolor=black!20!black,%
       linewidth=2pt,roundcorner=5pt,topline=true,
       frametitleaboveskip=\dimexpr-\ht\strutbox\relax,}
   \begin{mdframed}[]\relax%
}
{\end{mdframed}}

\newenvironment{cor}[1][]{%
\refstepcounter{teorem}%
\ifstrempty{#1}%
 {\mdfsetup{%
   frametitle={%
    \tikz[baseline=(current bounding box.east),outer sep=0pt]
    \node[anchor=east,rectangle,rounded corners,drop shadow,fill=black!40!brown]
         {\strut \textcolor{white}{Posljedica~\theteorem}};}}
 }
{\mdfsetup{%
  frametitle={%
   \tikz[baseline=(current bounding box.east),outer sep=0pt]
   \node[anchor=east,rectangle,rounded corners,drop shadow,fill=black!40!brown]
        {\strut \textcolor{white}{Posljedica~\theteorem:~#1}};}}%
 }%
\mdfsetup{innertopmargin=10pt,linecolor=black!20!brown,%
       linewidth=2pt,roundcorner=5pt,topline=true,
       frametitleaboveskip=\dimexpr-\ht\strutbox\relax,}
   \begin{mdframed}[]\relax%
}
{\end{mdframed}}

\newenvironment{primjer}
{\begin{Primjer}}{\hfill {\scalebox{0.6}{\stechak 5}} \end{Primjer}\vspace{-0.1cm}}

\newenvironment{vjezba}
{\begin{Vjezba}}{\end{Vjezba}\vspace{-0.1cm}}

\usepackage{imakeidx}
\makeindex

\renewcommand{\figurename}{Slika}
\addto\captionscroatian{
  
  }

\newcommand{\R}{\mathbb{R}}
\newcommand{\x}{\textbf{x}}
\newcommand{\sff}{{\textrm{I\!I}}}
\newcommand{\wtf}{{\textrm S}}

\newcommand{\fff}{{\textrm I}}
\newcommand{\n}{{\textbf{n}}}

\font\pleter=pleter
\font\stechak=stechak
\usepackage{picins}
\usepackage{microtype}
\usepackage{cutwin}

\newbox\mybox
\newdimen\myboxwidth
\newcommand\addpicture[3]{%
\setbox\mybox=\hbox{\includegraphics[scale=#3]{#2}}
\myboxwidth\wd\mybox
\renewcommand\windowpagestuff{%
\includegraphics[scale=#3]{#2}
}
\parpic[#1]{%
\begin{minipage}{\myboxwidth}
 \windowpagestuff
\end{minipage}
} }

\newcommand\addnumberedpicture[5]{%
\renewcommand{\figurename}{Sl.}
\setbox\mybox=\hbox{\includegraphics[scale=#3]{#2}}
\myboxwidth\wd\mybox
\renewcommand\windowpagestuff{%
\includegraphics[scale=#3]{#2}
\captionof{figure}{#4}\label{#5}
}
\parpic[#1]{%
\begin{minipage}{\myboxwidth}
 \windowpagestuff
\end{minipage}
}
\renewcommand{\figurename}{Slika}
}

\author{Vedad Pasic}
\title{Diferencijalna geometrija krivih i povr\v si u trodimenzionalnom euklidskom prostoru}
\date{Tuzla, 2017}

\newcommand\myrule{%
  \leavevmode\cleaders\hbox{\raisebox{-10pt}{\centerline{$\mbox{\stechak NNNNNNNNNNNNNNNNNNNNNN}$}}}\hfill\kern0pt}

  \newpagestyle{uvod}[\small]{
   
    \sethead[\tikz{
        \node[anchor=south,circle,draw,fill=black,text width=0.4cm,align=center]{\color{white}{\thepage}};}]
            [\chaptertitle]
            []
            {}
            {\chaptertitle}
            {\tikz{
        \node[circle,draw,fill=black,text width=0.4cm,align=center]{\color{white}{\thepage}};}}
   \setlength{\headsep}{1.5cm}

  }

  \newpagestyle{main2}[\small]{
    
    \sethead[\tikz{
        \node[anchor=south,circle,draw,fill=black,text width=0.4cm,align=center]
        {\color{white}{\thepage}};
        }
        ]
            [\chaptertitle]
            [
            \tikz{
            \node[rectangle,draw,fill=darkgray,text width=0.8cm,align=center]
            {\color{white}{\S \thesection}};
            }
            ]
            {
            \tikz{
        \node[rectangle,draw,fill=darkgray,text width=0.8cm,align=center]
        {\color{white}{\S \thesection}};
        }
            }
            {\sectiontitle}
            {
            \tikz{
            \node[circle,draw,fill=black,text width=0.4cm,align=center]
            {\color{white}{\thepage}};
            }
            }
   \setlength{\headsep}{1.5cm}
  }

\newpagestyle{title}[\small]{
    
    \sethead[\thepage]
            [{}]
            [{}]
            {}
            {}
            {\thepage}
       \setlength{\headsep}{1.5cm}
  }

  \newpagestyle{appendix}[\small]{
    
    \sethead[\tikz{
        \node[anchor=south,circle,draw,fill=black,text width=0.4cm,align=center]
        {\color{white}{\thepage}};
        }
        ]
            [\chaptertitle]
            [
            \tikz{
            \node[rectangle,draw,fill=darkgray,text width=0.8cm,align=center]
            {\color{white}{\thesection}};
            }
            ]
            {
            \tikz{
        \node[rectangle,draw,fill=darkgray,text width=0.8cm,align=center]
        {\color{white}{\thesection}};
        }
            }
            {\sectiontitle}
            {
            \tikz{
            \node[circle,draw,fill=black,text width=0.4cm,align=center]
            {\color{white}{\thepage}};
            }
            }
   \setlength{\headsep}{1.5cm}
  }

  \newpagestyle{appendix2}[\small]{
    
    \sethead[\tikz{
        \node[anchor=south,circle,draw,fill=black,text width=0.4cm,align=center]
        {\color{white}{\thepage}};
        }
        ]
            [\chaptertitle]
            [
            \tikz{
            \node[rectangle,draw,fill=darkgray,text width=0.8cm,align=center]
            {\color{white}{\thechapter}};
            }
            ]
            {
            \tikz{
        \node[rectangle,draw,fill=darkgray,text width=0.8cm,align=center]
        {\color{white}{\thechapter}};
        }
            }
            {\chaptertitle}
            {
            \tikz{
            \node[circle,draw,fill=black,text width=0.4cm,align=center]
            {\color{white}{\thepage}};
            }
            }
   \setlength{\headsep}{1.5cm}
  }

\makeatletter
\def\thickhrulefill{\leavevmode \leaders \hrule height 1ex \hfill \kern \z@}
\def\@makechapterhead#1{%
  \vspace*{5\p@}%
  {\parindent \z@ \centering \reset@font
        \thickhrulefill\quad
        \scshape \@chapapp{} \thechapter
        \quad \thickhrulefill
        \par\nobreak
        \vspace*{10\p@}%
        \interlinepenalty\@M
         {\scalebox{1}{{\pleter 5gggggggggggggg6}}}\\
        \Huge \bfseries #1\par\nobreak
        \par
         {\scalebox{1}{{\pleter 5gggggggggggggg6}}}\\
    \vspace*{25\p@}%
  }}
\def\@makeschapterhead#1{%
  \vspace*{5\p@}%
  {\parindent \z@ \centering \reset@font
        \par\nobreak
        \vspace*{10\p@}%
        \interlinepenalty\@M
        {\scalebox{1}{{\pleter 5gggggggggggggg6}}}\\
        \Huge \bfseries #1\par\nobreak
        \par
        {\scalebox{1}{{\pleter 5gggggggggggggg6}}}\\
    \vspace*{25\p@}%
  }}

\begin{document}
\WarningsOff*
\pagestyle{title}
\maketitle
\thispagestyle{empty}
\pagebreak
\pagenumbering{roman}
\section*{Differential Geometry of Curves and Surfaces in Three-dimensional Euclidean Space}
This manuscript, titled Differential Geometry of Curves and Surfaces in Three-dimensional Euclidean Space, is intended for undergraduate students of mathematics and other sciences with the need of use of Euclidean differential geometry. We deal with the differential geometry of curves and surfaces in three dimensions, covering all the main definitions, ideas and concepts needed for successful understanding and further application of this fundamental branch of mathematics. Many concrete examples with solutions are given with illustrations.
\section*{Diferencijalna geometrija krivih i povr\v si u trodimenzionalnom euklidskom prostoru}
Ova knjiga namijenjena je dodiplomskim studentima matematike i drugih nauka sa potrebom upotrebe euklidske diferencijalne geometrije. Upoznat ćemo se sa diferencijalnom geometrijom krivih i površi u tri dimenzije, pokrivajući sve glavne definicije, ideje i koncepte potrebne za razumijevanje i uspješnu dalju primjenu ove fundamentalne grane matematike. Dati su mnogi konkretni problemi sa rješenjima i mnoštvom ilustracija.

\vspace{2.5cm}

\begin{center}
{\scshape Sva prava zadr\v zana. Svako objavljivanje, \v stampanje ili umno\v zavanje
zahtjeva odobrenje autora}
\end{center}
\pagebreak
\section*{Diferencijalna geometrija na Univerzitetu u Tuzli}
{\bf Predmet}: Diferencijalna geometrija \\~ \\
{\bf Nastavnik}: Vedad Pa\v si\' c \\~ \\
{\bf Semestar}: VI (ljetni semestar III godine)  \\~ \\
{\bf ECTS}: 6 \\~ \\
{\bf Kabinet}: PMF 313 \\~ \\
{\bf Email}: vedad.pasic@untz.ba \\~ \\
{\bf Web}: http://www.pmf.untz.ba/vedad/nastava.html \\ ~ \\
{\bf Facebook}: Grupa "Diferencijalna geometrija - Univerzitet u Tuzli" \\~ \\ ~ \\~ \\

{\bf{Organizacija}}\\ ~
\begin{itemize}
\item 3h predavanja (srijeda 12-15) i 2h vje\v zbi (petak 10-12).
\item Sedmi\v cne problemske zada\' ce diskutovane na vje\v zbama...
\item ...rad na zada\' ci je predispitna obaveza i nosi 15 \% ukupne ocjene.
\item Predispitne obaveze - test 1 (krive  - 20\%), test 2 (povr\v si  - 20\%) i zada\' ce i aktivnost (15\%).
\item Finalni ispit - 45\%.
\end{itemize}
\dominitoc[n]
\begingroup
\makeatletter \let\ps@plain\ps@empty \makeatother
\tableofcontents
\endgroup

\newpage
\pagenumbering{arabic}
\pagestyle{uvod}

\makeatletter
\let\ps@plain\ps@empty
\makeatother

\chapter{Uvod}
\vspace{1cm}
Osnovna misija ove knjige i predmeta kojem je namijenjena je
prou\v cavanje geometrije krivih i povr\v si koriste\' ci se metodama diferencijalnog ra\v cuna (primarno izvodima i integralima), odakle i dolazi naziv ``{\it diferencijalna\/} geometrija''. U na\v sem izu\v cavanju geometrije krivih i povr\v si, bit \' cemo zainteresovani samo za one osobine istih koje su nezavisne od njihove pozicije u prostoru, tojest za one osobine koje su invarijantne pod pomjeranjima euklidskog $3$-prostora, kao \v sto su translacije i rotacije. U tu svrhu jako mnogo direktne primjene u diferencijalnog geometriji ima i linearna algebra.

Ova knjiga treba da \v citaoca pripremi za ozbiljnije teme iz diferencijalne geometrije, uvode\' ci sve osnovne koncepte ove matemati\v cke oblasti za krive i povr\v si u trodimenzionalnom euklidskom prostoru.

Ve\' c du\v zi niz godina vodim se idejom da treba napisati modernu knjigu iz diferencijalne geometrije na maternjem jeziku, jer takvog ud\v zbenika jednostavno na \v zalost nema. Nije realisti\v cno o\v cekivati od studenta dodiplomskog studija da izu\v cava ovu tematiku koriste\' ci se pregr\v stom razli\v citih njiga, posebno ne na stranim jezicima.

Ova knjiga je poku\v saj da se taj problem premosti i da studenti dodiplomskog studija matematike, ali i drugih studija kojima je ova tematika interesantna, mogu diferencijalnu geometriju izu\v cavati na kvalitetetan, moderan i ugodan na\v cin.

Rije\v c `\emph{geometrija}' dolazi od gr\v ckog
$\gamma\epsilon\omega\mu\epsilon\tau\rho\acute\iota\alpha$
i slo\v zena je rije\v c od ``$\gamma\tilde{\eta}$'' = ``geo'' - Zemlja i
``$\mu\epsilon\tau\rho\acute\omega\nu$'' = ``metron'' - mjeriti.
Stoga je ``geometrija'' potekla od nauke mjerenja zemlje - naprimjer iz potrebe
mjerenja povr\v sine poljoprivrednih polja.

Najstariji zabilje\v zeni po\v ceci geometrije se mogu na\' ci u drevnoj Mesopotamiji i Egiptu u drugom milenijumu prije nove ere. Rana geometrija je bila skup empiri\v cki otkrivenih principa koji su se ticali du\v zina, uglova, povr\v sina i zapremina, koje su se razvile kako bi mogle biti primijenjene u geodetskim mjerenjima, gradnji, astronomiji i raznim zanatima.

Najraniji poznati tesktovi o geometriji su egipatski Rhindov matemati\v cki papirus (2000-1800 pne) i Moskovski papirus (1890 pne), te babilonske glinene tablice, kao \v sto je Plimpton 322 tablica (1900 pne). Naprimjer, Moskovski papirus daje formulu za izra\v cunavanje zapremine zarubljene piramide (\emph{frustuma}), jer ko o \v cemu, Egip\' cani o piramidama.

Stari Grci su ovu `primijenjenu nauku' pretvorili  u teorijski matemati\v cku granu, prou\v cavaju\' ci geometrijske objekte u apstraktnom smislu (Euklidovi ``Elementi'' su najpoznatiji primjer prvog ud\v zbenika iz ove oblasti, a koristio se u u\v cionicama sve do pro\v slog vijeka).

\emph{Diferencijalna} geometrija  razvila se mnogo kasnije, tojest nakon \v sto su Newton i Leibnitz razvili diferencijalni i integralni ra\v cun, u ranom $18$-om vijeku: tada je bilo mogu\' ce prou\v cavati kompliciranije geometrijske objekte, kao \v sto na primjer su proizvoljno zakrivljene
``krive'' i ``povr\v si'' u $3$-prostoru.

Dosta ovog predmeta \' ce biti posve\' ceno opisivanju ``zakrivljenosti''
objekata (krivih i povr\v si) u prostoru, ili generalnije, njihovom
``obliku''.

O `krivoj' se mo\v ze misliti kao o obliku koji bi poprimila savijena \v zica
u prostoru (ili na ravni);  podjednako, mo\v zemo o njoj misliti kao
o tragu koji ostavlja \v cestica koja se kre\' ce u prostoru.
Jednostavni i ve\' c vi\dj{}eni primjeri su prava linija
ili kru\v znica u prostoru (ili na ravni).

O ``povr\v si'' se mo\v ze misliti kao o sapunici ili mjehuri\' cu sapuna
ili povr\v sini tijela na Zemlji. Jednostavan primjer je ravan ili sfera
u prostoru.

Ove primitivne pojmove ``krive'' i ``povr\v si'' \' cemo preciznije definisati kad za to do\dj{}e vrijeme.
\begin{primjer}
Posmatrajmo elipsu:
$\displaystyle
   E = \left\{(x,y) \in \mathbb{R}^2\,\left|\,\left(\frac{x}{a}\right)^2 + \left(\frac{y}{b}\right)^2 = 1\right.\right\}$.  Elipsa je data u \emph{implicitnoj formi}, tojest pomo\' cu jedna\v cine.
Mo\v zemo li je napisati u \emph{parametarskoj formi},  tojest u formi
$$E=\{\gamma(t)\,|\, t\in I\},$$ gdje je $I$ neki interval?
 Prije svega, mo\v zemo je predstaviti {\it kao graf:\/}
   $E\supset\{(x,\pm b\,\sqrt{1-(\frac{x}{a})^2})\,|\,-a<x<a\}$
   --- ali na ovaj na\v cin ne dobijemo cijelu elipsu
   (rije\v site za $y$ da biste dobili preostale ta\v cke;
    primijetite da ne mo\v zemo imati $x=\pm a$ bez gubitka
    diferencijabilnosti).

{\it U parametarskoj formi:\/} elipsu je data sa
   $\displaystyle E=\{(a\cos t,b\sin t)\,|\,t\in\R\}.$
   Me\dj{}utim, na ovaj na\v cin pokrivamo elipsu beskona\v cno mnogo puta.
   Primijetite da $\gamma'\neq0$ svugdje, zbog osobina trigonometrijskih funkcija.
\end{primjer}
\begin{remark}
Teorema implicitnog preslikavanja (vidi Dodatak, teorem \ref{TeoremImplicitnogPreslikavanja}) garantira da svaka
kriva ili povr\v s data u implicitnoj formi mo\v ze, pod odre\dj{}enim
pretpostavkama, biti napisana u parametarskoj formi (i obratno!).
Teorema implicitnog preslikavanja daje nam samo teorijski alat -- obi\v cno nije prakti\v cno da je se koristi kako bi se stvarno  {\it na\v sla\/}  parametrizacija krive ili povr\v si.
Me\dj{}utim, ona nam daje vrijedan kriterij za kada (par) jedna\v cina opisuje
povr\v s (krivu): to bi trebao biti nivo skup {\it submerzije\/}
$\R^3\to\R$ (tojest~$\R^2\to\R$).
\end{remark}
\begin{remark}
Nivo skup diferencijabilne funkcije $f:\R^n\to\R$ koji
odgovara vrijednosti $c\in \mathbb{R}$ je skup ta\v caka $\{x\in \R^n\ |\ f(x)=c\}$.
\end{remark}
\begin{primjer}
Nivo skup
funkcije $f(x,y,z)=x^2+y^2+z^2$ vrijednosti $c$ je sfera $x^2+y^2+z^2=c$ sa centrom
u $(0,0,0)$ i radijusom $\sqrt{c}$.
\end{primjer}

Diferencijalna geometrija ima niz primjena u prirodnim naukama, a
najzna\v cajnije su primjene u fizici, naprimjer kada posmatramo krivu
kao trag koji ostavlja \v cestica koja se kre\' ce kroz prostor.

Diferencijalna geometrija je jezik Einsteinove generalne teorije relativnosti i cijeli svemir je iskazan i obja\v snjen konceptima i metodama diferencijalne geometrije. Razumijevanje koncepta zakrivljenosti je klju\v cno u izra\v cunavanju pozicija satelita u orbiti, a diferencijalne geometrija se koristi i u izu\v cavanju gravitacionih so\v civa i crnih rupa. Diferencijalne forme, koncept izrastao iz ove nauke, koriste se kako bi se opisao elektromagnetizam, a imamo primjene u Lagrangeovoj mehanici i Hamiltonovoj mehanici.
Geometrijsko modeliranje, \v sto uklju\v cuje ra\v cunarsku grafiku, te CAD (computer--aided geometric design) se zasnivaju na idejama iz diferencijalne geometrije.
Niz drugih primjena u in\v zenjerstvu (prilikom rje\v savanja problema u obradi digitalnih signala), ekonomiji (ekonometrija, vidi naprimjer \cite{pasic2016matematika}), ra\v cunarskoj nauci, teoriji kontrole, vjerovatno\' ci, itd.

Ovaj ud\v zbenik je namijenjen studentima tre\' ce godine dodiplomskog studija matematike i treba Vas pripremiti za materiju koja \' ce se, izme\dj{}u ostalog, raditi na predmetima \v cetvrte godine Vi\v sa geometrija \cite{pasic2016visa} i Matemati\v cke metode u fizici \cite{pasic2016matematicke}.
Diferencijalna geometrija je u centru mog istra\v ziva\v ckog rada - a iskreno se nadam da \' ce neki od vas postati moje istra\v ziva\v cke kolegice i kolege. Ako vas to interesuje, pogledajte moje radove iz ove i drugih oblasti \cite{barakovic2017physical,karasuljic2017construction,okicic2002acting,pasic2009new,pasic2010new,pasic2014pp,
pasic2015torsion,pasic2017axial,pasic2014new,pasic2015uniformly,pasic2002differential,
pasic2005pp,rekic2015continuity}.

Potrebna predznanja:
\begin{itemize}
\item Osnovna linearna algebra i analiti\v cka geometrija;
\item (Vektorski) diferencijalni i integralni ra\v cun.
\end{itemize}
\subsubsection{Notacija}
\begin{itemize}
\item $u\cdot v$ - za skalarni proizvod dva vektora;
\item $||u||=\sqrt{u\cdot u}$ - za du\v zinu (normu, intenzitet) vektora;
\item $u\times v$ - za vektorski proizvod dva vektora;
\item $d_pf$ - za Jacobijevu matricu funkcije
   $f$ u ta\v cki $p$;
\item $f_u$, $f_v$ - za parcijalne izvode funkcije $f=f(u,v)$;
\item[{\scalebox{0.6}{\stechak L}}] Zavr\v setak dokaza.
\item[{\scalebox{0.6}{\stechak K}}] Zavr\v setak primjedbe.
\item[{\scalebox{0.6}{\stechak 5}}] Zavr\v setak primjera.
\item[{\scalebox{0.6}{\stechak O}}] Zavr\v setak rje\v senja.
\end{itemize}
\subsubsection{Literatura}
Kako bi se kompletirala ova knjiga, kori\v steno je jako mnogo literature, ali najva\v zniji
ud\v zbenici su oni klasi\v cni, kao \v sto je Struik  \cite{struik2012lectures}
ili Spivak \cite{spivak1999comprehensive} u kojim se mo\v ze na\' ci velika ve\' cina
materijala koji pokriva ova knjiga. Dodatno, za specifi\v cnije stvari, jako
su korisni i ud\v zbenici do Carma \cite{do1976differential}, Kuhnela
\cite{kuhnel2015differential}, Valirona \cite{valiron1986classical}
i Sharipova \cite{sharipov2004course}, te rad Sya \cite{sy2001general} o op\' cim heliksima. Na jezicima ju\v znih Slavena ima jako
malo literature o diferencijalnoj gemetriji koja je ozbiljna i/ili komprehensivna,
a mogu se pogledati naprimjer Stojanovi\' c  \cite{stojanovic1963osnovi},
Mihajlovi\'c \cite{mihajlovic1968elementi}, \v Zarinac-Fran\v cula \cite{zarinac1990diferencijalna}, Bla\v zi\' c i Bokan \cite{blazic1998uvod}, dok je dobar resurs za diferencijalnu geometriju krivih i povr\v si u Wolfram Mathematica Velimirovi\' c \cite{velimirovic2010geometrija} kao i Abbena i drugi \cite{abbena2017modern}. Vi\v se o historiji matematike, a samim tim umonogome i geometrije, mo\v zete nau\v citi od Radoj\' ci\' ca \cite{radojcic1950opsta}. Za predznanja iz matemati\v cke analize toplo preporu\v cujem Dedagi\'ca \cite{dedagic2005matematicka}
ili bilo koju drugu kvalitetnu knjigu iz maemati\v cke analize.

\subsubsection{Sadr\v zaj ukratko}
\begin{itemize}
\item Krive (na ravni i u $3$-prostoru);
\item Elementarna teorija povr\v si;
\item Elementarne jedna\v cine povr\v si;
\item Geometrija na povr\v si.
\end{itemize}

~ \newpage ~ \hspace{1cm} ~  \par ~ \\

\vfill {\scalebox{0.7}{{\stechak lkkkkkkkkkkkkkkkkkd}}}

\chapter{Diferencijalna geometrija prostornih krivih}
\pagestyle{main2}
\minitoc
\vspace{1cm}
U ovom poglavlju uvest \' cemo koncept parametrizovane, geometrijske i implicitne \emph{krive} u trodimenzionalnom euklidskom prostoru, te pridru\v zene geometrijske objekte, du\v zinu luka, tangentni, normalni i binormalni vektor, te \' cemo diskutovati o geometrijskim veli\v cinama svake krive iz kojih mo\v zemo povratiti originalnu krivu - krivinu i torziji.

\section{Analiti\v cka reprezentacija}
\begin{defn}[Kriva]\label{defn:kriva} \index{kriva!definicija}
\emph{Kriva} je slika jednog glatkog (ili mogu\' ce nekoliko glatkih) preslikavanja $$\gamma:I\to\R^3,$$ gdje je $I$ realni interval.
Kriva je \emph{regularna} ako $\gamma'(t)\neq0,$ $\forall\  t\in I$.
\end{defn}
\begin{remark}
Pod glatkim preslikavanjem dakako podrazumjevamo preslikavanje klase $C^k$ za neko $k\in \mathbb{N}$. Nadalje pretpostavljamo da su sve krive $$\gamma:I\to\R^3$$ najmanje klase $C^3$ ako se ne naglasi druga\v cije.
\end{remark}
\begin{remark}
Uvijek \' cemo pretpostavljati da su na\v se krive regularne, ukoliko nije navedeno druga\v cije. \index{kriva!regularna}
Ovaj dogovor osigurava da sve na\v se krive imaju dobro definisane tangentne linije: \index{tangentne linije}
\begin{equation} \label{eqn:tangente}
   u \mapsto \gamma(t) + u\,\gamma'(t).
\end{equation}
\end{remark}
\begin{defn}[Koordinate krive]\label{defn:koordinateKrive}
Komponente $x$, $y$ i $z$ ta\v cke $(x,y,z)\in\R^3$
\' cemo zvati njenim \emph{koordinatama}.\newline
Kada je $$t\mapsto\gamma(t)=(x(t),y(t),z(t))$$ (regularna) kriva,
onda \' cemo funkcije $x(t)$, $y(t)$, i $z(t)$ zvati
\emph{koordinatne funkcije}.
\end{defn}
\renewcommand{\figurename}{Sl.}
\begin{primjer}[Prava linija]\index{prava}
Svaka prava $t\mapsto\gamma(t)$ mo\v ze biti parametrizovana kao
$$x(t)=a_1+tb_1, \
      y(t)=a_2+tb_2, \
      z(t)=a_3+tb_3, \ t\in\R,$$
\addnumberedpicture{r}{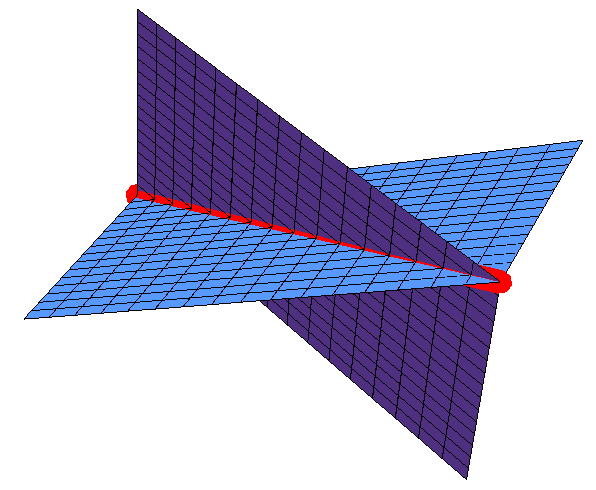}{0.25, bb= 0 0 596 490}{Prava presjek ravni}{fig:pravaRavni}
    gdje su $a_i,b_i$ konstante i barem jedno $b_i\neq0$. Ova kriva je stoga jasno regularna, jer je stoga $\gamma'(t) = (b_1,b_2,b_3) \ne 0$ po pretpostavci.
   Pretpostavimo, bez gubitka op\' cenitosti, da je $b_1\neq0$; onda mo\v zemo eliminisati parametar $t$ i napisati krivu kao graf preko cijele $x$-ose:
\begin{eqnarray*}
      y(x)&=&a_2+\frac{b_2}{b_1}(x-a_1), \\
      z(x)&=&a_3+\frac{b_3}{b_1}(x-a_1),\ \ x\in\R
\end{eqnarray*} \v Sta predstavljaju ove jedna\v cine?
\end{primjer}

\begin{primjer}[Kru\v znica]\label{primjer:kruznica}\index{kru\v znica} Kru\v znica u prostoru, tojest u ravni $z=0$ centrirana u koordinatnom po\v cetku $t\mapsto \gamma(t)$ je data sa :
\addnumberedpicture{l}{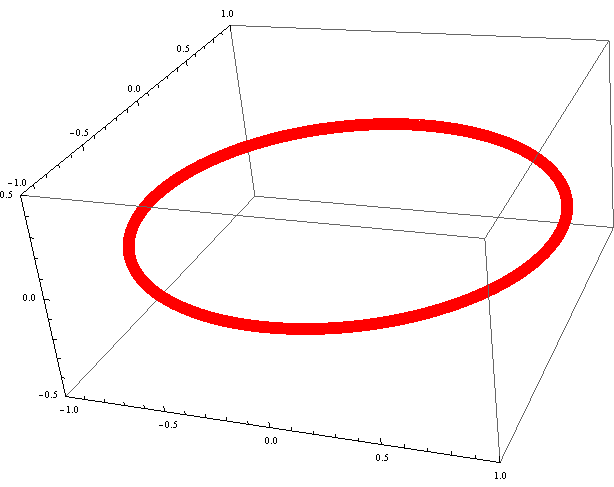}{0.2, bb=0 0 614 487}{Kru\v znica}{fig:kruznica}
\begin{eqnarray*}
      x&=&r\cos t, \\
      y&=&r\sin t, \\
      z&=&0, \qquad t\in[0,2Pi].
\end{eqnarray*}
Ova kriva je \emph{planarna} krrva,\index{kriva!planarna} tojest nalazi se u ravni. Ovo je o\v cito regularna kriva, jer je $$\gamma'(t) = r(-\sin t, \cos t, 0) \ne 0, \ \forall t.$$ Stoga, ako primjenimo jedna\v cinu (\ref{eqn:tangente}) u ovom slu\v caju, imamo da je tangentna linija u proizvoljnoj ta\v cki $\xi\in\R$ na kru\v znicu data sa
\[
t\mapsto r(\cos \xi - t\sin \xi, \sin \xi + t \cos \xi, 0).
\] Stoga je tangenta na jedini\v cnu kru\v znicu u $xOy$ ravni centriranoj u koordinatnom po\v cetku u ta\v cki $t=\frac\pi4$, tojest $ \left(\frac{1}{\sqrt{2}}, \frac{1}{\sqrt{2}}, 0 \right)$, data sa
\[
t\mapsto \frac{1}{\sqrt{2}} (1-t,1+t,0).
\] Vidi Sliku \ref{fig:tangente} za vi\v se grafi\v ckih primjera.
\begin{figure}[h!]
\centering
\includegraphics[height=4cm,bb = 0 0 987 848]{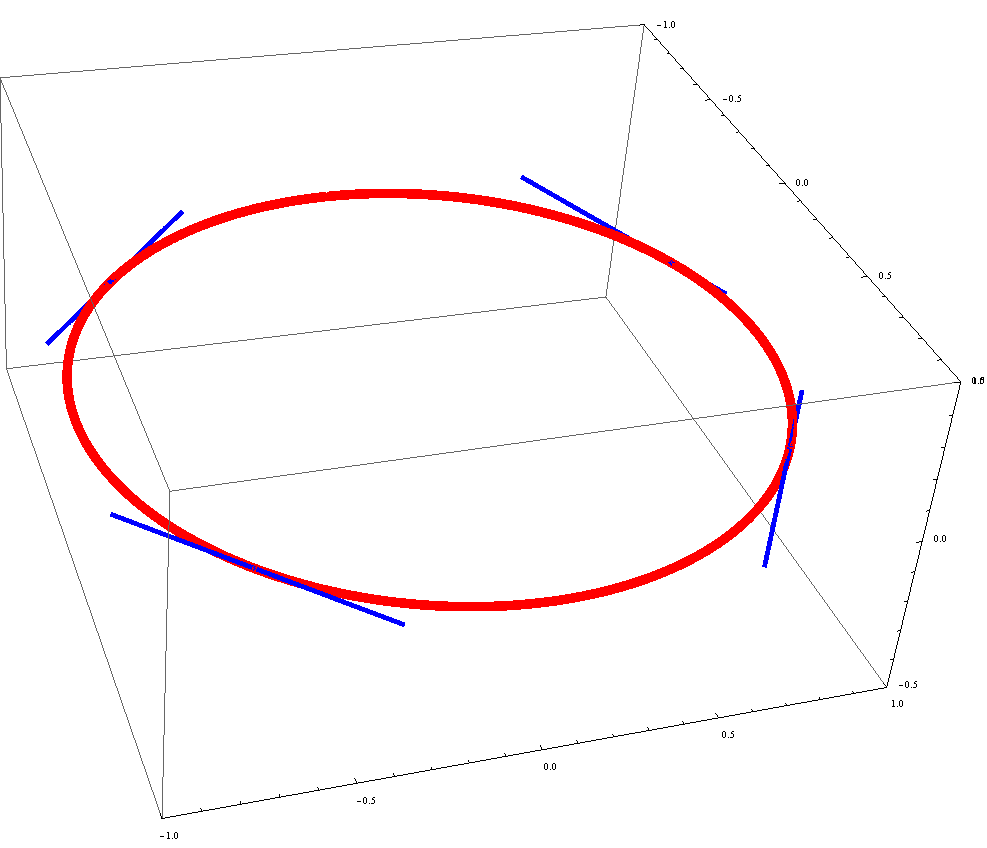}
\caption{Tangente na jedini\v cnu kru\v znicu u ta\v ckama $t=1,2,4,5.5$}\label{fig:tangente}
\end{figure}
\end{primjer}
\begin{remark}
Mogu\' ce je da dvije razli\v cite krive imaju isti kodomen. Ako se vratimo primjeru \ref{primjer:kruznica} i posmatrajmo novu krivu $t\mapsto \gamma_1(t) :$
\[
\gamma_1(t) = (\cos t, \sin t, 0), \qquad t\in[0,4Pi].
\]  Slika ove krive je ponovo o\v cito jedini\v cna kru\v znica koja le\v zi u ravni $z=0$ centrirana u koordinatnom po\v cetku, ali nije injekcija za razliku od krive iz primjera \ref{primjer:kruznica}.
\end{remark}

\addnumberedpicture{r}{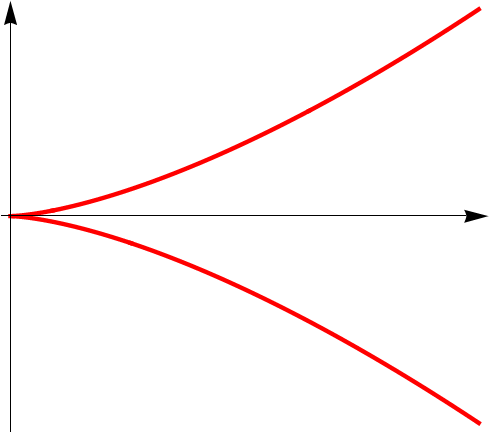}{0.25,bb=0 0 489 432}{Iregularna kriva}{fig:iregularna}
Ako kriva nije regularna, onda se naziva \emph{iregularnom}.
Takva je naprimjer
\[
t\mapsto \gamma(t), \qquad \gamma(t) = (t^2,t^3,0), t\in \R.
\] O\v cito je da je $$\gamma'(t) = (2t,3t^2,0)$$ i stoga je $\gamma'(t)$ jednako nula vektoru za $t=0$. Stoga kriva $\gamma$ nije regularna, vidi Sliku \ref{fig:iregularna}.

\begin{primjer}[Kru\v zni heliks] Ova vrlo posebna kriva se defini\v se pomo\' cu jedna\v cina \begin{eqnarray*}
      x&=&r\cos t, \\
      y&=&r\sin t, \\
      z&=&h\ t,
\end{eqnarray*}
\addnumberedpicture{r}{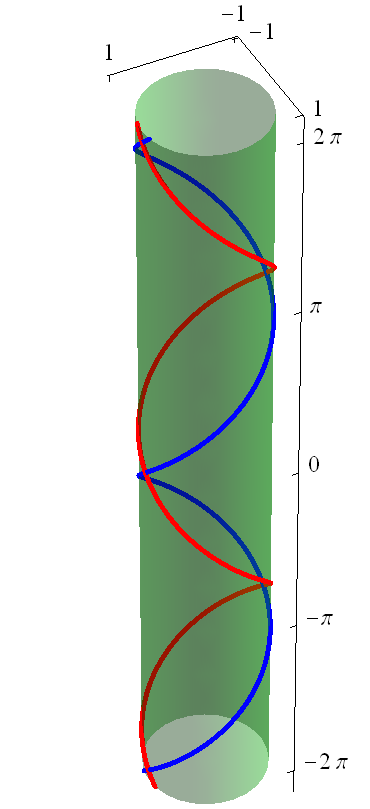}{0.2,bb = 0 0 370 809}{Heliks}{fig:lijeviDesniHeliks}
gdje su $r$ i $h$ realni brojevi razli\v citi od nule, a $t \in \R$ parametar.

Kao \v sto \' cemo vidjeti, ova kriva \' ce se iznova pojavljivati u primjerima i primjenama, te \' cemo posvetiti cijelu podsekciju klasi helikoidnih krivih, vidi sekciju \ref{sek:heliksi}.

Ova kriva le\v zi na cilindru $x^2+y^2=r^2$, te se na istom \emph{uvija} na na\v cin da kada se parametar $t$ pove\' ca za $2\pi$, koordinate $x$ i $y$ se vrate na istu vrijednost, dok se koordinata $z$ pove\' ca za $2\pi h$, \v sto nam daje vertikalnu udaljenost heliksovih zavoja, tojest \emph{nagib} heliksa.

Kada je konstanta $h$ pozitivna, tada imamo \emph{desni} kru\v zni heliks (crvena kriva na slici \ref{fig:lijeviDesniHeliks}), a kada je konstanta $h$ negativna, tada imamo \emph{lijevi} kru\v zni heliks (plava kriva na slici \ref{fig:lijeviDesniHeliks}).
\addnumberedpicture{l}{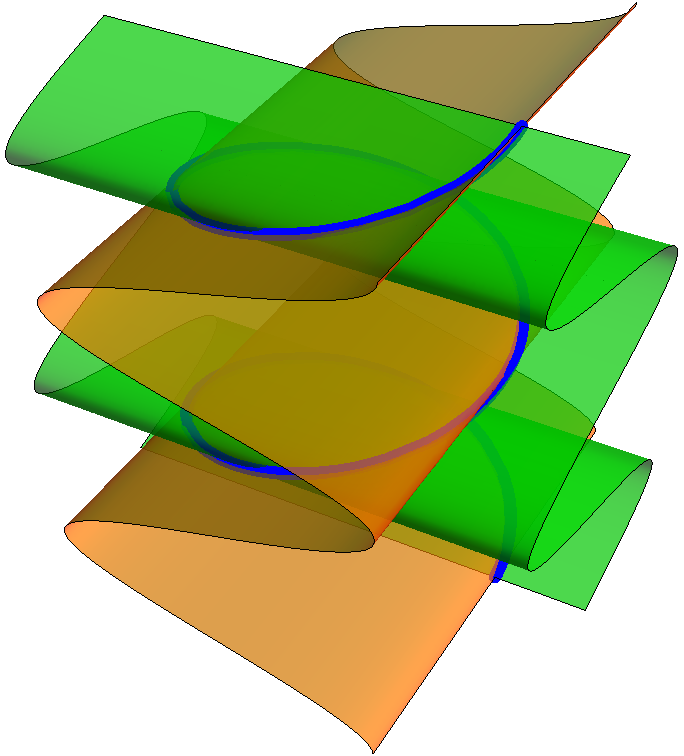}{0.2,bb=0 0 684 756}{Heliks implicitno}{fig:presjekHeliks}
Ovakva podjela helika je nezavisna od izbora parametara ili parametrizacija - ona je intrinzi\v cna osobina kru\v znih heliksa. Lijevi heliks se nikad ne mo\v ze preslikati na desni i obratno, \v sto zna svako ko je ikad radio sa vijcima ili konopcima.

Ako je $h\neq0$, onda heliks mo\v zemo napisati kao graf preko njegove ose ($z$-ose):
\begin{eqnarray*}
      x&=&r\cos\frac{z}{h}, \\
      y&=&r\sin\frac{z}{h}, z\in\R,
\end{eqnarray*}
vidi sliku \ref{fig:presjekHeliks} kako presjek povr\v si $x=\cos z$ i $y=\sin z$ daju kru\v zni heliks radijusa $1$ i nagiba $2\pi$.
\end{primjer}%

\begin{defn}[Reparametrizacija krive]\label{defn:reparametrizacijaKrive}
\emph{Reparametrizacija} krive $(a,b)\ni t\mapsto\gamma(t)\in\R^3$
je nova parametrizirana kriva
\begin{equation}
\tilde\gamma (s)=\gamma(u(s)),
\end{equation}
gdje je $u:(\tilde a,\tilde b)\mapsto (a,b)$ sirjekcija i $u'\ne0$.
\end{defn}
Uslov $u'(s)\neq0$ $\forall \ s$ osigurava da
su reparametrizacije regularnih krivih regularne, \v sto proizilazi iz lan\v canog pravila
\begin{equation}\label{eqn:regularnostReparametrizacije}
\tilde{\gamma}'(s) = u'(s) \cdot \gamma'(u(s)) \ne 0.
\end{equation}
Preslikavanje $u:(\tilde a,\tilde b)\mapsto (a,b)$  iz Definicije \ref{defn:reparametrizacijaKrive} je \emph{difeomorfizam}, vidi Definiciju \ref{defn:imerzija}.
\begin{defn}[Ekvivalentnost krivih]\label{defn:ekvivalentneKrive}
Za dvije krive $\alpha$ i $\beta$ ka\v zemo da su \emph{ekvivalentne} i pi\v semo $\alpha\equiv \beta$ ukoliko postoji difeomorfizam $u:(\tilde a,\tilde b)\mapsto (a,b)$ klase $C^k$ takav da je $\beta = \alpha \circ u$.
\end{defn}
Uslov regularnosti (\ref{eqn:regularnostReparametrizacije}) implicira da je $u'$ stalnog znaka, a za reparametrizaciju ka\v zemo da \v cuva orijentaciju, odnosno da je \emph{pozitivna reparametrizacija} ukoliko je $u'>0$, a da je \emph{negativna reparametrizacija} ukoliko je $u'<0$.
\begin{primjer}
Reparametrizirajmo krivu \[
t\mapsto \gamma(t), \quad \gamma(t) := (e^t \cos 2t, 2, e^t \sin 2t), \ t \in (-\infty,\infty).
\] Posmatrajmo funkciju
\[
u: (0,\infty) \mapsto (-\infty,\infty), \quad u(s) = \ln \left(\frac{s}{\sqrt{5}}\right).
\]
Zamjenom $t=u(x)$ u originalnoj parametrizaciji, dobijamo
\[
\tilde{\gamma} (s)  = \left(
\frac{s}{\sqrt{5}} \cos \ln \frac{s^2}{5}, 2, \frac{s}{\sqrt{5}} \sin \ln \frac{s^2}{5}
\right)
\]  Ovo se naravno mo\v ze \v ciniti nepotrebnim, ali ima velike prednosti, kako \' cemo vrlo brzo vidjeti!
\end{primjer}
\begin{defn}[Geometrijska kriva]\label{defn:geometrijskaKriva}
Klasa ekvivalencije $[\gamma]$ je \emph{geometrijska} (neparametrizovana) kriva.
\end{defn}
\begin{remark}
Definicija \ref{defn:geometrijskaKriva} nam omogu\' cava distinkciju izme\dj{}u \emph{krive} i \emph{parametrizacije krive}. Svaka (geometrijska) kriva mo\v ze imati beskona\v cno mnogo parametrizacija, razli\v citih svojstava. Lokalno posmatrano, krive su ekvivalentne ukoliko imaju iste skupove slika.
\end{remark}

\subsection{Implicitna forma krive}
\emph{U punoj verziji teksta.}

\subsection*{Vje\v zbe}
\begin{vjezba}
Na\' ci tangentu na kru\v zni heliks $t\mapsto (\cos t,\sin t, t)$ u ta\v cki $t=\frac{\pi}{3}$ i napisati je u kanonskom obliku.
\end{vjezba}
\begin{vjezba}\label{vjezba:implicitnaKriva}
Posmatrajmo krivu datu implicitno sa
$$
4x^2+9y^2+36z^2=36,\qquad \sqrt{\frac{3}{5}}\,z = \frac{x}{3},
$$
Provjerite regularnost krive sa slike \ref{vjezba:implicitnaKriva} i parametrizirajte je.
\end{vjezba}
\begin{figure}
\centering
\includegraphics[width=0.7\textwidth, bb =0 0 676 421]{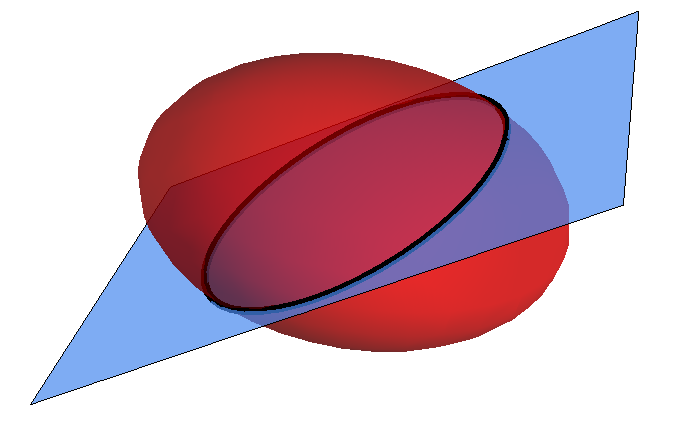}
\caption{Implicitno zadata kriva iz vje\v zbe \ref{vjezba:implicitnaKriva}}
\end{figure}
\begin{vjezba}
Poka\v zite da su koni\v cni presjeci regularne krive i na\dj{}ite njihove parametrizacije.
\begin{enumerate}
\item[(a)] $x^2+y^2=z^2$, $x+z=\sqrt{2}$;
\item[(b)] $x^2+y^2=z^2$, $x+\sqrt{3}z=2$;
\item[(c)] $x^2+y^2=z^2$, $\sqrt{3}x+z=2$.
\end{enumerate}
\end{vjezba}
\begin{vjezba}
Doka\v zite teoremu implicitnog preslikavanja pomo\' cu teoreme inverznog preslikavanja, koriste\' ci se na\v som notacijom.
\end{vjezba}
\begin{vjezba}\label{vjezba:Picard}
Doka\v zite Picard-Lindel\" ofovu teoremu: Neka su $I$ otvoreni potskup $\mathbb{R}$ i $U$ otvoreni potskup $\mathbb{R}^n$ i neka je $I\times U\ni(x,y)\mapsto f(x,y)\in\mathbb{R}^n$ neprekidna i Lipschitz neprekidna po $y$ i $(x_0,y_0)\in I\times U$;
onda {\it problem po\v cetne vrijednosti\/}
$$
   y'(x) = f(x,y(x)), \quad y(x_0)=y_0
$$
ima jedinstveno lokalno rje\v senje ako je $f$ diferencijabilna funkcija.
\end{vjezba}

\section{Du\v zina luka krive}
\emph{U punoj verziji teksta.}

\subsection*{Vje\v zbe}
\begin{vjezba}\label{vjezba:invarijantnostLuka}
Pokazati da je du\v zina luka invarijantna pod reparametrizacijom krive.
\end{vjezba}
\begin{vjezba}
Izra\v cunajte du\v zinu luka i reparametrizirajte, ako je mogu\' ce, du\v zinom luka krive $t \mapsto \gamma_i$ :
\begin{enumerate}
\item[(a)]    $\displaystyle \gamma_1(t) = \left(
   \frac{3\sqrt{3}}{2\sqrt{2}} \cos(t), 2\sin(t), \frac{\sqrt{5}}{2\sqrt{2}}\cos(t)
   \right)$
\item[(b)] $\displaystyle \gamma_2(t)= \frac{\sqrt{2}}{2}\left(1-\frac{t^2}{2}, \sqrt{2}t, 1+\frac{t^2}{2}\right)$.
\item[(c)] $\displaystyle \gamma_3(t)= \left(
1-\cosh(t),
\sinh(t),
t
\right)$.
\item[(d)] $\gamma_4(t) = e^t(\cos t,\sin t,0)$
\item[(e)] $\displaystyle \gamma_5(t) = (e^t \cos 2t, 2, e^t \sin 2t)$
\end{enumerate}
\end{vjezba}

\section{Stripovi i oskulatorna ravan}
\emph{U punoj verziji teksta.}

\subsection*{Vje\v zbe}
\begin{vjezba}\label{vjezbe:linearnaZavisnost}
Potvrdite detaljno da je $t\mapsto\gamma(t)$ prava linija ako su
$\gamma''(t)$ i $\gamma'(t)$ linearno zavisne za sva $t$.
\end{vjezba}
\begin{vjezba}\label{vjezbeOskulatornaHeliks}
Izra\v cunati oskulatornu ravan na kru\v zni heliks u proizvoljnoj ta\v cki.
\end{vjezba}
\begin{vjezba}[Bonus]

\begin{enumerate}
\item[(a)] Iskoristiti Wolfram Mathematica, kako biste napravili funkciju koja ra\v cuna jedna\v cinu oskulatorne ravni \\ $paramOskRavan[\gamma\_, t\_, t0\_]$ na krivoj $\gamma$ parametra $t$ u ta\v cki $t0$ u parametarskom obliku.
\item[(b)] Napraviti animaciju oskulatorne ravni na kru\v znom heliksu, kao \v sto je to ura\dj{}eno na slici na predmetnoj stranici.
\end{enumerate}
\end{vjezba}

\section{Krivina i torzija}
\emph{U punoj verziji teksta.}

\subsection{Krivina}

\emph{U punoj verziji teksta.}

\subsection{Frenet--Serret principalni okvir}
\emph{U punoj verziji teksta.}

\subsection{Torzija}
\emph{U punoj verziji teksta.}

\subsection*{Vje\v zbe}
\begin{vjezba}
Po definiciji izra\v cunati krivinu i torziju krivih $t\mapsto \gamma_i (t)$:
\begin{itemize}
\item $\gamma_1= (e^t\cos t, e^t\sin t,0)$;
\item $\gamma_2= (t,\cosh t,1)$;
\item $\gamma_3= (\cos^2 t,\sin^2 t,\sin^2 t +1)$;
\item $\gamma_4= \left(\frac45 \cos t, 1+\sin t, 2-\frac35\cos t\right)$;
\end{itemize}
\end{vjezba}
\begin{vjezba}\label{vjezbe:oskulatornaKruznica}
Uvjerite se da se oskulatorna kru\v znica du\v zinom luka parametrizovane krive $\gamma$ u ta\v cki $\gamma(s_0)$ (koristite $s_0=0$ kako bi vam bilo lak\v se) doista mo\v ze parametrizovati (du\v zinom luka) tako da se Taylorovi polinomi drugog reda zaista poklapaju.
\end{vjezba}
\begin{vjezba}\label{vjezbe:formulaKrivinaTorzija}
Neka je $t\mapsto\gamma(t)$ parametrizovana kriva tako da $T'(t)\neq0,$
$\forall t$. Pokazati da mo\v zemo izabrati $N=\frac{T'}{|T'|}$ i da je, sa ovakvim izborom
$N$,
\begin{equation*}
   \kappa = \frac{|\gamma'\times\gamma''|}{|\gamma'|^3},\qquad
   \tau = \frac{|\gamma',\gamma'',\gamma'''|}{|\gamma'\times\gamma''|^2}.
\end{equation*}
\end{vjezba}
\begin{vjezba}
Izra\v cunati krivinu i torziju krive $t\mapsto\gamma(t)$ definisane sa $\displaystyle  \gamma(t)= \left(t,\frac{1+t}{t},\frac{1-t^2}{t}\right)$ koriste\' ci se formulama. Dajte geometrijski argument za\v sto $\tau$ uzima vrijednost koju uzima. Poku\v sajte izra\v cunati krivinu po definiciji.
\end{vjezba}

\begin{figure}[h]
\centering
\includegraphics[width=\textwidth,bb=0 0 1293 595]{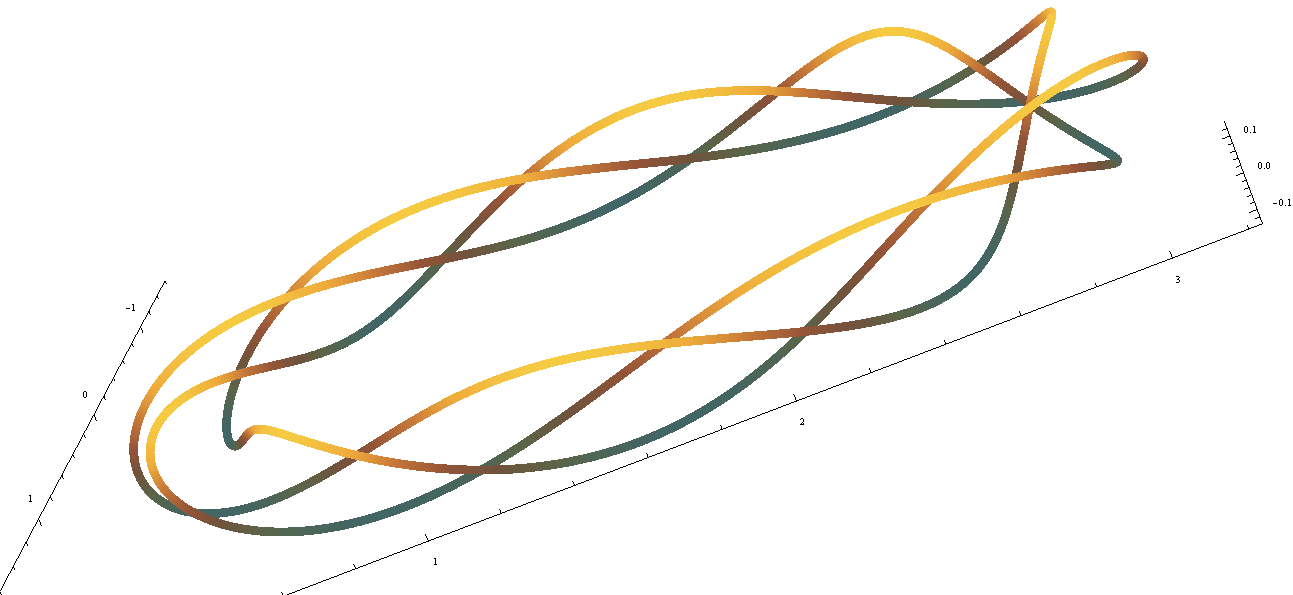}
\caption{Kriva iz vje\v zbe \ref{vjezba:TorzijaJedan}}\label{fig:torzijajedan}
\end{figure}

\begin{vjezba}\label{vjezba:TorzijaJedan}
Izra\v cunati krivinu i torziju krive $u\mapsto \gamma(u)$ \v cije su parametarske funkcije date sa
\begin{align*}
x(u) &= -\sqrt{2}\sin\frac{3 u}{2} + \frac{2\sqrt{2}}{21} \sin\frac{7 u}{2} - \frac{\sqrt{2}}{39} \sin\frac{13 u}{2}, \\
y(u)&= \frac{102\sqrt{2}}{91} -\sqrt{2} \cos\frac{3 u}{2} -\frac{2\sqrt{2}}{21}  \cos\frac{7 u}{2} - \frac{\sqrt{2}}{39} \cos\frac{13 u}{2}, \\
z(u) &= \frac{2}{15} \sin 5u.
\end{align*} \v Sta je to posebno u ovom primjeru, datom na slici \ref{fig:torzijajedan}?
\end{vjezba}

\begin{vjezba}
Neka je sada $t\mapsto \gamma(t)$ regularna parametrizovana kriva. Poka\v zite da krivina zadovoljava $\kappa^2=|\gamma'\times\gamma''|^2$, kada je $N={T'\over{|T'|}}$.

Neka je sada $s\mapsto \tilde{\gamma}(s)=\gamma(t(s))$
reparametrizacija krive $\gamma$ (ne obavezno du\v zinom luka). Poka\v zite da je, za sva $s$,
$$
{{|\tilde{\gamma}'\times \tilde{\gamma}''|^2}\over{|\tilde{\gamma}'|^6}}(s)
={{|{\gamma}'\times {\gamma}''|^2}\over{|{\gamma}'|^6}}(s)
$$
Zaklju\v cite da je, za regularnu parametriziranu krivu $t\mapsto \gamma(t)$,
$$
\kappa^2={|\gamma'\times \gamma''|^2\over |\gamma'|^6}.
$$
\end{vjezba}

\begin{vjezba}[Bonus]
\begin{enumerate}
\item[(a)] Napraviti funkcije $krivina[g\_,t\_]$ i \\ $torzija[g\_,t\_]$ u Wolfram Mathematici koje ra\v cunaju krivinu i torziju date proizvoljno parametrizovane krive $g$ parametra $t$.
\item[(b)] Napravite animaciju trobrida $(T,N,B)$ na proizvoljnoj krivoj. Opcionalno, prika\v zite i vrijednosti krivine i torzije prilikom animacije.
\end{enumerate}
\end{vjezba}


\section{Frenet-Serret jedna\v cine}
\emph{U punoj verziji teksta.}

\subsection*{Vje\v zbe}
\begin{vjezba}
Popunite `praznine' u dokazu Frenetovih jedna\v cina te poka\v zite da doista jeste
\begin{equation*}
   F' = F\,
   \left(\begin{array}{ccc}
      0 & -\kappa & 0 \\
      \kappa & 0 & -\tau \\
      0 & \tau & 0\end{array}\right)
   \quad\Leftrightarrow\quad
   \left\{\begin{array}{cccc}
      T'=& & \kappa N& \\
      N'=&-\kappa T& &+\tau B\\
      B'=& &-\tau N&
   \end{array}\right..
\end{equation*}
\end{vjezba}
\begin{vjezba}
Neka $F=(T,N,B)$ ozna\v cava principalni okvir du\v zinom luka parametrizovane krive
$s\mapsto\gamma(s)$ i neka su $\kappa(s)$ i $\tau(s)$ krivina i torzija
krive $\gamma$.
Defini\v simo ``Darboux-ovo vektorsko polje'' $$
   \delta := N\times N' = \tau(s)\,T + \kappa(s)\,B.
$$
Pokazati da su Frenetove jedna\v cine ekvivalentne jedna\v cinama
\begin{equation*}
   T' = \delta\times T, \quad N' = \delta\times N, \quad  B' = \delta\times B.
\end{equation*}
\end{vjezba}
\begin{vjezba}
Pokazati da kriva $s\mapsto\gamma(s)$ uzima vrijednosti na sferi ako i samo ako
njene normalne ravni prolaze kroz (fiksnu) ta\v cku $c\in\mathbb{R}^3$ (centar sfere).
\end{vjezba}
\begin{vjezba}
\begin{enumerate}
\item[(a)] Neka je $t\mapsto\beta(t)\in S^2\subset\mathbb{R}^3$ regularna kriva.

    Defini\v simo krivu
\begin{equation}
   t\mapsto\gamma(t):=\int_{t_0}^t\beta(t)\times\beta'(t)dt.\tag{\dag}
\end{equation}
Pokazati da $\gamma$ ima konstantnu torziju.
\item[(b)] Sada pretpostavimo da kriva $s\mapsto\gamma(s)$ (parametrizovana du\v zinom luka)
ima konstantnu torziju $\tau\equiv1$.
Pokazati da je $\gamma$ oblika $(\dag)$ sa odgovaraju\' com krivom
$s\mapsto\beta(s)\in S^2$.
\end{enumerate}
\end{vjezba}
\begin{vjezba}\label{vjezba:FTPK}
Neka je $s\mapsto\kappa(s)$ data realna funkcija. Defini\v simo funkciju
$\varphi(s):=\int_{s_0}^s\kappa(s)ds$.
Provjerite da $$
   \gamma(s) :
   = \left(\int_{s_0}^{s}\cos\varphi(s)ds,\int_{s_0}^{s}\sin\varphi(s)ds,0\right)
$$
defini\v se du\v zinom luka parametrizovanu planarnu krivu sa krivinom $\kappa(s)$.
\end{vjezba}

\subsection{Heliksi}\label{sek:heliksi}
\emph{U punoj verziji teksta.}

\subsubsection{Generisanje op\' ceg heliksa}
\emph{U punoj verziji teksta.}

\subsection*{Vje\v zbe}
\begin{figure}
\centering
\includegraphics[height=6cm,bb = 0 0 602 660]{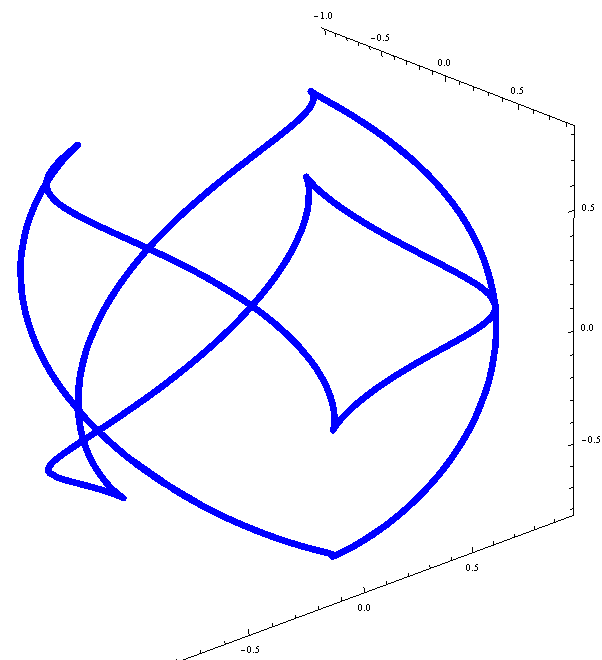}
\caption{Op\' ci heliks iz vjezbe \ref{vjezba:Heliks} }\label{fig:opciheliks3}
\end{figure}

\begin{vjezba}
Da li je kru\v zni heliks op\' ci heliks? Za\v sto?
\end{vjezba}
\begin{vjezba}
Poka\v zite da je klasa krivih $t\mapsto \gamma(t)$ data sa
$$
\gamma(t) = \left(\frac{bt^n}{n}, \sqrt{2b} \frac{t^{n+1}}{n+1}, \frac{t^{n+2}}{n+2} \right),
$$ gdje je $b>0$ realna konstanta, a $n\in \mathbb{N}$, klasa op\' cih heliksa.
\end{vjezba}
\begin{vjezba}\label{vjezba:Heliks}
Doka\v zite da je kriva $t\mapsto \gamma(t)$, vidi sliku \ref{fig:opciheliks3} data sa
$$
\gamma(t)=\left(\frac45 \cos 2t - \frac15 \cos 8t, \frac45 \sin 2t - \frac15 \sin 8 t,  \frac \cos 3t\right),
$$ gdje je $\pi/3<t<2\pi/3$, op\' ci heliks.
\end{vjezba}
\begin{vjezba}[Bonus]
Kreirajte u Wolfram Mathematica funkciju koja provjerava da li je data kriva heliks, te funkciju koja iz date planarne krive generi\v se op\' ci heliks.
\end{vjezba}

\section{Prirodne jedna\v cine i fundamentalna teorema prostornih krivih}
\emph{U punoj verziji teksta.}

\section{Okvirne krive}
\emph{U punoj verziji teksta.}

\subsection*{Vje\v zbe}
\begin{vjezba}
Izvesti dokaz strukturnih jedna\v cina za okvirne krive.
\end{vjezba}
\begin{vjezba}
Izvesti dokaz Fundamentalnog teorema za okvirne krive.
\end{vjezba}
\begin{vjezba}\label{vjezba:zaParalelni}
Posmatrajmo krivu $t\mapsto \gamma(t)$ datu sa $\gamma(t)=(t\cos t,t\sin t,t)$.
\begin{enumerate}
\item[(a)] Doka\v zite da je kriva $\gamma$ regularna i da je
$$
t\mapsto N(t) = {1\over{\sqrt{2}}} (-\cos t,-\sin t,1)
$$ jedini\v cno normalno vektorsko polje za $\gamma$.
\item[(b)] Na\dj{}ite adaptirani okvir $F=(T,N,B)$ za strip $(\gamma,N)$ i izra\v cunajte njegovu torziju $\tau$, geodezijsku krivinu $\kappa_g$ i normalnu krivinu $\kappa_n$.
\end{enumerate}
\end{vjezba}

\begin{vjezba}
Neka $\gamma$ parametrizira pravu, tojest neka je $\gamma'\times \gamma''\equiv 0$. Neka je $F$ bilo koji adaptirani okvir krive $\gamma$.

Poka\v zite da su $\kappa_n = \kappa_g = 0$. Na\dj{}ite strip $(\gamma,N)$ takav da je $\kappa_n = \kappa_g = 0$ i  $\tau \equiv 1$.
\end{vjezba}
\begin{vjezba}\label{vjezba:rotacijaokvira}
Neka je $F=(T,N,B)$ adaptirani okvir i neka je $\tilde F=FA$, gdje je $$
   A = \left(\begin{array}{ccc}
      1 & 0 & 0 \\
      0 & \cos\varphi &-\sin\varphi \\
      0 & \sin\varphi & \cos\varphi \end{array}\right)
$$
sa nekom funkcijom $t\mapsto\varphi(t)$.
\begin{itemize}
\item Uvjerite se da je $\tilde F$ jo\v s jedan adaptirani okvir.
\item Izra\v cunajte kako se strukturne jedna\v cine (tojest $\kappa_n$,
$\kappa_g$ i $\tau$) mijenaju.
\end{itemize}
\end{vjezba}

\subsection{Paralelna normalna polja i paralelni okviri}
\emph{U punoj verziji teksta.}

\subsection*{Vje\v zbe}
\begin{vjezba}
Na\dj{}ite paralelno jednini\v cno normalno vektorsko polje $\tilde{N}$ za krivu $\gamma$ iz vje\v be \ref{vjezba:zaParalelni}. Mo\v zete ga zapisati u obliku
$$
\tilde{N}(t) = \cos(\varphi(t)) N(t) + \sin(\varphi(t))B(t),
$$ gdje je $\varphi(t)$ funkcija po $t$ koju trebate odrediti.

[Pomo\' c: $\int \frac{1}{\sqrt{t^2+2}} dt = arcsinh \frac{t}{\sqrt{2}} + C$]

Potom izra\v cunati $\kappa_g$ i $\kappa_n$ za paralelni strip krivine $(\gamma,\tilde{N})$.
\end{vjezba}

\begin{vjezba}
Doka\v zite da su bilo koja dva paralelna adaptirana okvira krive $t\mapsto\gamma(t)$
povezana konstantnom rotacijom u normalnoj ravni.
\end{vjezba}
\begin{vjezba}
Poka\v zite da kriva uzima vrijednosti na sferi ako i samo ako krivine
$\kappa_g$ i $\kappa_n$ paralelnog okvira zadovoljavaju jedna\v cinu
planarne krive. \par
Kako se mo\v ze pro\v citati radijus sfere iz ove jedna\v cine?
\end{vjezba}

\subsection{Frenetovi okviri}
\emph{U punoj verziji teksta.}

\subsection*{Vje\v zbe}
\begin{vjezba}
 Neka je kriva $t\mapsto\alpha(t)$ data sa
$$
\alpha(t) = \big( {t+\sin t}, {-t+\sin t}, \sqrt{2}{\cos t}\big).
$$
\begin{enumerate}
\item[(a)] Na\dj{}ite reparametrizaciju du\v zinom luka krive $s\mapsto\alpha(s)$.
\item[(b)] Na\dj{}ite Frenetov okvir $F=(T,N,B)$ krive $\alpha$.
\item[(c)] Izra\v cunajte krivinu i torziju krive $\alpha$ po definiciji.
\item[(d)] \v Sta mo\v zete zaklju\v citi o krivoj $\alpha$ iz njene krivine i torzije i na osnovu \v cega?
    \end{enumerate}
\end{vjezba}

\vfill {\scalebox{0.7}{{\stechak lkkkkkkkkkkkkkkkkkd}}}
\chapter{Diferencijalna geometrija povr\v si}
\pagestyle{main2}
\minitoc
\vspace{1cm}
U ovom poglavlju uvest \' cemo koncept parametrizovane i implicitne \emph{povr\v si} u trodimenzionalnom euklidskom prostoru, te pridru\v zene geometrijske objekte, prvu i drugu fundamentalnu formu, Gaussovo preslikavanje, te opisati oblik povr\v si pomo\' cu razli\v citih definicija krivina. Razmatrat \' cemo i razli\v cite posebne klase povr\v si.
\section{Analiti\v cka reprezentacija povr\v si}
\begin{defn}[Povr\v s]
\emph{Povr\v s} je slika jednog (ili mogu\' ce nekoliko) glatkog
preslikavanja $$\x:U\to\R^3,$$ gdje je $U\subset\R^2$. Povr\v s je \emph{regularna}
ako je $\x$ imerzija, tojest ako $d_{(u,v)}\x:\R^2\to\R^3$ je injekcija
za sva  $(u,v)\in U$, vidi Definiciju \ref{defn:imerzija}.
\end{defn}
\begin{remark}
Povr\v s $\x$ je regularna ako i samo ako su $\x_u(u,v)$ i
$\x_v(u,v)$ linearno nezavisni za sva $(u,v)\in U$, tojest
ako i samo ako $
   \x_u\times\x_v \neq 0.
$
\end{remark}
\underbar{\bf Dogovor.}
Uvijek \' cemo pretpostaviti da su na\v se povr\v si regularne.
\begin{primjer}[Ravan]
   Tri nekolinearne ta\v cke $P_0, P_1, P_2 \in \R^3$, nalaze se
   \addnumberedpicture{r}{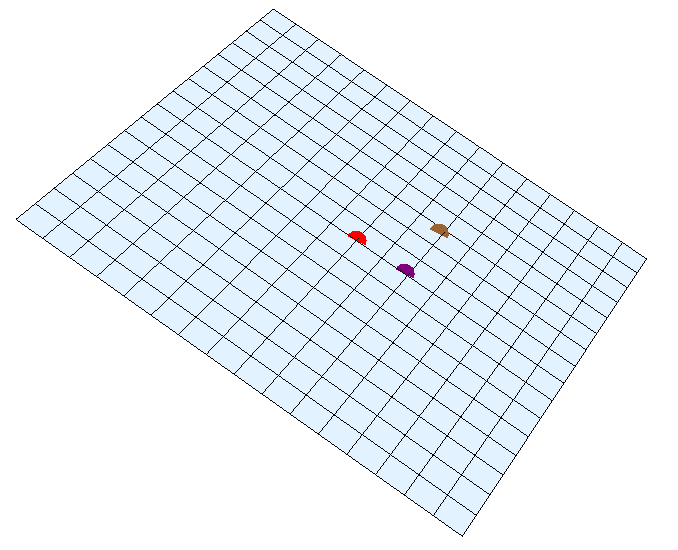}{0.2, bb=0 0 678 550}{Ravan}{fig:ravan}
   na jedinstvenoj ravni, parametrizovanoj sa
      $$(u,v)\mapsto P_0+u\,(P_1-P_0)+v\,(P_2-P_0),$$
vidi Dodatak \ref{appendix:predznanje}. Povr\v s je regularna jer je
\[
x_u \times x_v = (P_1-P_0) \times (P_2-P_0) \ne 0,
\]  jer su ta\v cke $P_0,P_1,P_2\in\R^3$ nekolinearne.
\end{primjer}
\begin{primjer}[Sfera]
\v Cesta parametrizacija je data sa
\addnumberedpicture{r}{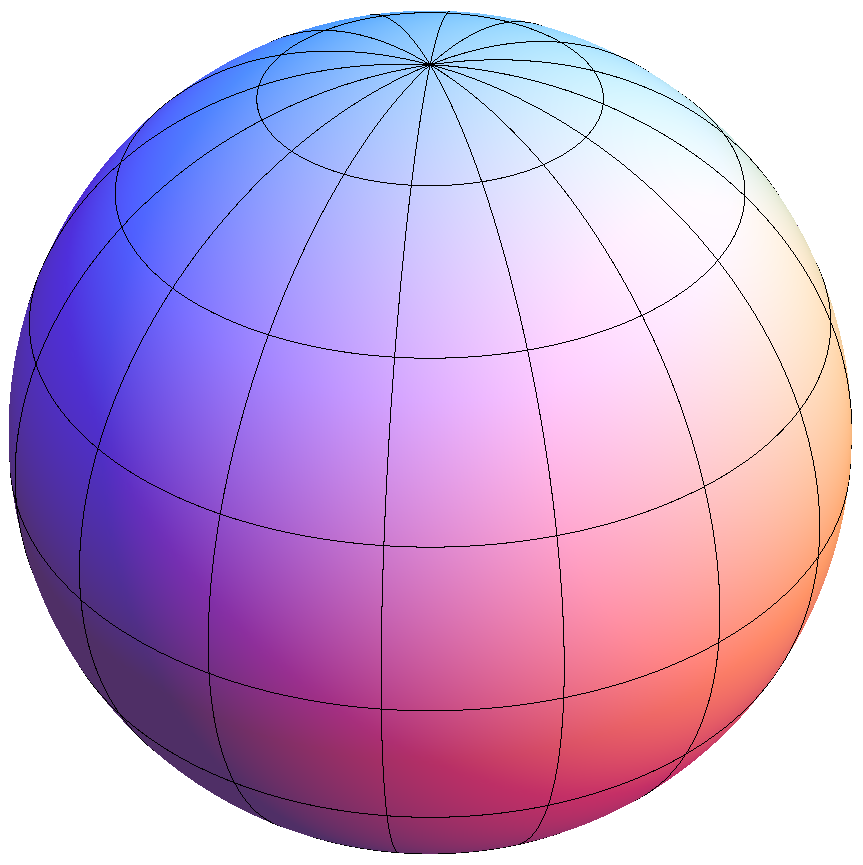}{0.15, bb= 0 0 863 865}{Sfera}{fig:sfera}
      $$(u,v)\mapsto(\cos u\cos v,\cos u\sin v,\sin u),$$
$u \in \left[-\frac{\pi}{2},\frac{\pi}{2}\right], v \in [0,2\pi]$. Me\dj{}utim, imamo problem \v sto se ti\v ce regularnosti.

Naime, parametrizacija prestaje da bude regularna  za
$\cos u=0$ i $\sin u=\pm1$ (``sjeverni'' i ``ju\v zni polovi'' sfere); ovaj problem je simptomati\v can i ne mo\v ze biti rije\v sen:
{\it ne postoji regularna parametrizacija cijele sfere odjednom\/}.

Stoga, treba nam nekoliko parametrizacija kako bismo pokazali da je to regularna povr\v s
(na primjer, uzmite gornju parametrizaciju i sli\v cnu parametrizaciju koja je rotirana tako da
njeni ``polovi'' le\v ze na $(x,y)$-ravni):
\[
(u,v)\mapsto(\cos u\sin v,\sin u\sin v,\cos v),
\] $u \in \left[0,\pi\right], v \in [0,2\pi]$.
\end{primjer}
\subsubsection{Potsjetnik na sferne koordinate}
Sferni koordinatni sistem ima za ideju orijentaciju na sfernoj povr\v sini. Sferu mo\v zemo podijeliti "paralelnim"
kru\v znicama, me\dj{}u kojima je i ekvatorijalna, koje nazivamo paralelama i "velikim" kru\v znicama koje sve prolaze
kroz polove sfere, koje nazivamo meridijanima.

U takvoj podjeli sfere, pokazuje se boljim opis polo\v zaja ta\v cke u
smislu koliko smo daleko (u stepenima) od nekog fiksnog meridijana i koliko smo daleko (u stepenima) od neke fiksne
paralele, od uobi\v cajenih koordinata, du\v zine, \v sirine i visine.

Naravno, tre\' ca
bitna stvar o polo\v zaju je i udaljenost od koordinatnog po\v cetka. Postoje dva pristupa sfernim koordinatama, u
zavisnosti koju paralelu biramo za fiksnu.

\begin{figure}[!h]
\centering
\includegraphics[width=0.8\textwidth,bb = 0 0 1019 557]{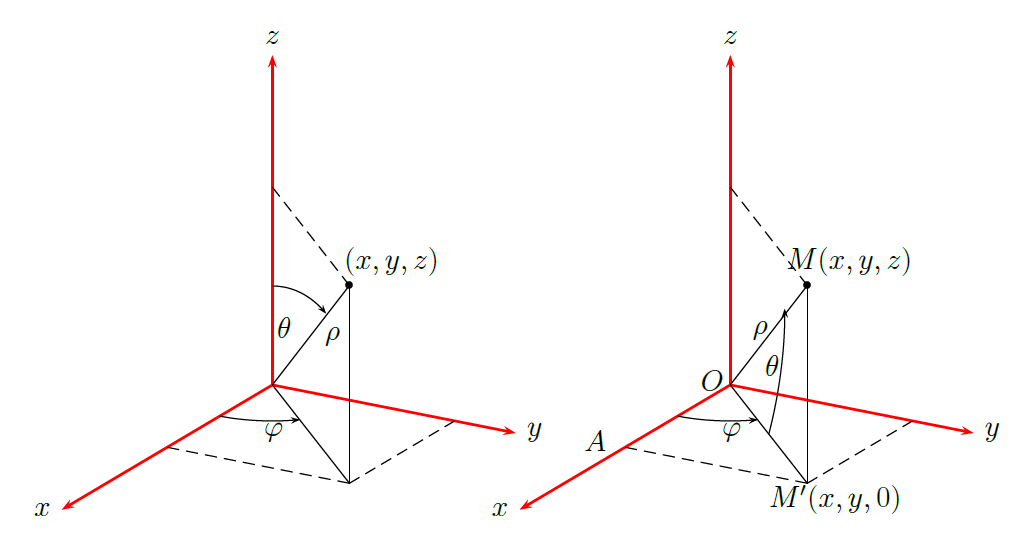}
\caption{Dva pristupa sfernim koordinatama}\label{slikasfernekoordinate}
\end{figure}

Posmatrajmo sliku \ref{slikasfernekoordinate} desno. Uzimamo da je $\rho$ udaljenost ta\v cke $M$ od koordinatnog po\v cetka, $\varphi$ je
udaljenost od meridijana, tojest ugao izme\dj{}u potega $OM'$, gdje je $M'$ projekcija ta\v cke $M$ u $Oxy$ ravan i
pozitivnog dijela $x$-ose i $\theta$ je ugao izme\dj{}u $\rho$ i $Oxy$ ravni, tojest udaljenost ta\v cke $M$ od ekvatorijalne ravni.

Uo\v cimo trougao $\triangle OM'M$. To je pravougli trougao, pa iz njega o\v citavamo
$$\cos \theta =\frac{OM'}{OM} \ , \ \sin \theta =\frac{MM'}{OM} \ ,$$
odnosno
\begin{equation}\label{sfernekoordinate1}
OM'=\rho \cos \theta \ , \ MM'=z=\rho \sin \theta \ .
\end{equation}
Iz pravouglog trougla $\triangle OAM'$ o\v citavamo
$$\cos \varphi =\frac{OA}{OM'} \ , \ \sin \varphi =\frac{AM'}{OM'} \ ,$$
tojest
\begin{equation}\label{sfernekoordinate2}
OM'=\frac{x}{\cos \varphi} \ , \ OM'=\frac{y}{\sin \varphi} \ .
\end{equation}
Kombinuju\' ci (\ref{sfernekoordinate1}) i (\ref{sfernekoordinate2}) dobijamo sferne koordinate za slu\v caj kada ugao
$\theta$ mjerimo od ekvatorijalne ravni.
$$x=\rho \cos \varphi \cos \theta \ , \ y=\rho \sin \varphi \cos \theta \ , \ z=\rho \sin \theta \ .$$
Prirodne granice su
$$0\leq \rho <+\infty \ , \ 0\leq \varphi \leq 2 \pi \ , \ -\frac{\pi}{2}\leq \theta \leq \frac{\pi}{2} \ ,$$
jer su od "ekvatora" najudaljeniji polovi i to sjeverni $90^\circ$, a ju\v zni $-90^\circ$. 

Uvode\' ci sferni koordinatni sistem, mjere\' ci ugao $\theta$ od sjevernog pola (Slika \ref{slikasfernekoordinate},
lijevo), sistem glasi
$$x=\rho \cos \varphi \sin \theta \ , \ y=\rho \sin \varphi \sin \theta \ , \ z=\rho \cos \theta \ ,$$
gdje su prirodne granice novih koordinata
$$0\leq \rho <+\infty \ , \ 0\leq \varphi \leq 2 \pi \ , \ 0\leq \theta \leq \pi \ .$$
Sada je najudaljenija ta\v cka od sjevernog pola, ju\v zni pol i to $180^\circ$. 

\begin{primjer}[(Kru\v zni) Helikoid]\label{primjer:helikoid}
Ovo je regularna povr\v s, iscrtana
\setbox\mybox=\hbox{\includegraphics[scale=0.15, bb = 0 0 640 1050]{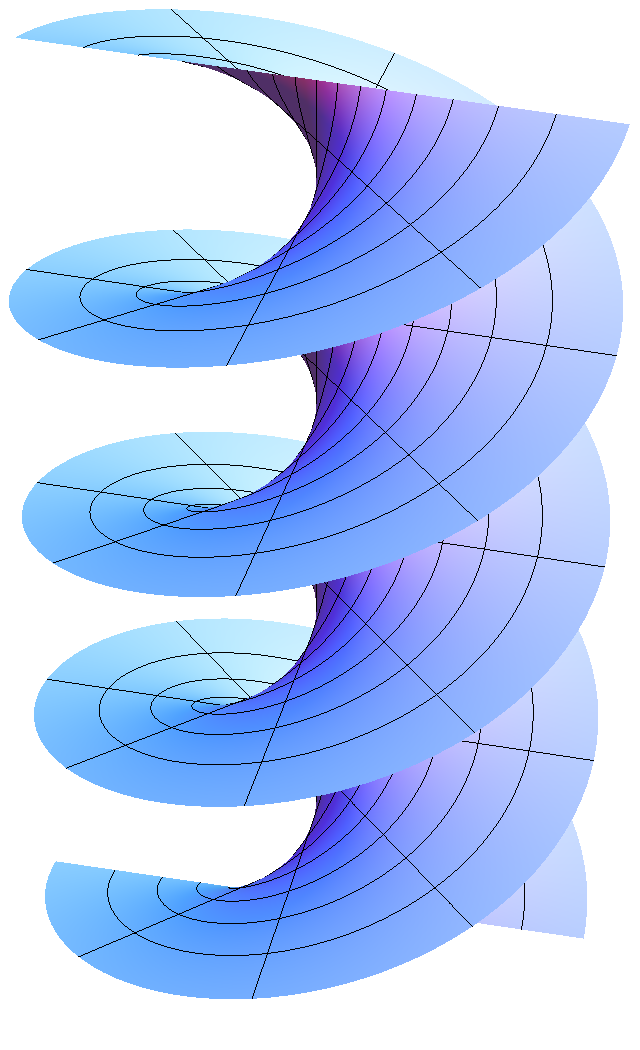}}
\myboxwidth\wd\mybox
\renewcommand\windowpagestuff{%
\includegraphics[scale=0.15, bb = 0 0 640 1050]{helikoid.png}
\captionof{figure}{Helikoid}\label{fig:helikoid}
\vspace{0.5cm}
}
\parpic[l]{%
\begin{minipage}{\myboxwidth}
 \windowpagestuff
\end{minipage}
}

pravom normalnom linijom koja se kre\' ce du\v z heliksa:
      $$(u,v)\mapsto\x(u,v):=(\sinh u\cos v,\sinh u\sin v,c\ v), \ c\in \mathbb{R}$$
Povr\v s je regularna jer je za sva $(u,v)\in\R^2$ $$
      (\x_u\times\x_v)(u,v)
      = \cosh u\,(\sin v,-\cos v,\sinh u)
      \neq 0.
   $$
   Ovo je veoma posebna povr\v s.
Kru\v zni helikoid je minimalna povr\v s koja ima kru\v zni heliks kao granicu. Ona je jedina linijska minimalna povr\v s osim ravni, a \v cak do 1992. godine je bio jedini znani primjer kompletne minimalne embedirane povr\v si kona\v cne topologije sa beskona\v cnom krivinom, dok nije otkrivena Hoffmannova minimalna povr\v s.
\end{primjer}
\begin{defn}
Kad je povr\v s data sa $\x=(x,y,z)$ onda se funkcije $x(u,v)$, $y(u,v)$, $z(u,v)$
zovu \emph{koordinatne funkcije} povr\v si.
\end{defn}
\begin{defn}
\emph{Reparametrizacija} povr\v si $U\ni(u,v)\mapsto\x(u,v)\in\R^3$
je nova parametrizovana povr\v s$$
   \tilde\x(\tilde u,\tilde v) = \x(\phi(\tilde u,\tilde v)),
$$
gdje je $\phi:\tilde U\to U$ difeomorfizam, vidi Definiciju \ref{defn:imerzija}.
\end{defn}
\begin{primjer}
Kru\v zni helikoid, slika \ref{fig:helikoid},  iz primjera \ref{primjer:helikoid} o\v cito se mo\v ze jednostavnije reparametrizirati sa
$$(u,v)\mapsto\x(u,v):=(u\cos v,u\sin v,c\ v), \ c\in \mathbb{R}.$$
\end{primjer}
\begin{remark}
Ako je $\tilde\x(\tilde u,\tilde v) = \x(\phi(\tilde u,\tilde v))$, gdje je
$\phi(\tilde u,\tilde v) = (u(\tilde u,\tilde v),v(\tilde u,\tilde v))$,
onda je
$$
   \tilde\x_{\tilde u}\times\tilde\x_{\tilde v}
   = \left| \begin{array}{cc}
   u_{\tilde u} & u_{\tilde v} \\
   v_{\tilde u} & v_{\tilde v}
   \end{array}
   \right|\,
      (\x_u\times\x_v)\circ\phi.
$$
Stoga je reparametrizacija regularne povr\v si regularna.
\end{remark}
\subsection{Implicitna forma povr\v si}
Posmatrajmo skup$$
   \Sigma := \{(x,y,z)\,|\,F(x,y,z)=0\}.
$$
Pretpostavimo da je ${\partial F\over\partial z}(x,y,z)\neq0$ za sva $(x,y,z)\in C$.
Onda (po teoremi implicitnog preslikavanja \ref{TeoremImplicitnogPreslikavanja}) mo\v zemo {\it lokalno\/} rije\v siti jedna\v cinu za $z=z(x,y)$, tojest $$
   C\cap B = \{(x,y,g(x,y))\,|\,(x,y)\in U\}
$$
za neke okoline $B\subset\R^3$ od $(x,y,z)$ i $U\subset\R^2$.

\begin{figure}[!h]
\centering
\includegraphics[height=4cm,bb = 0 0 1320 1649]{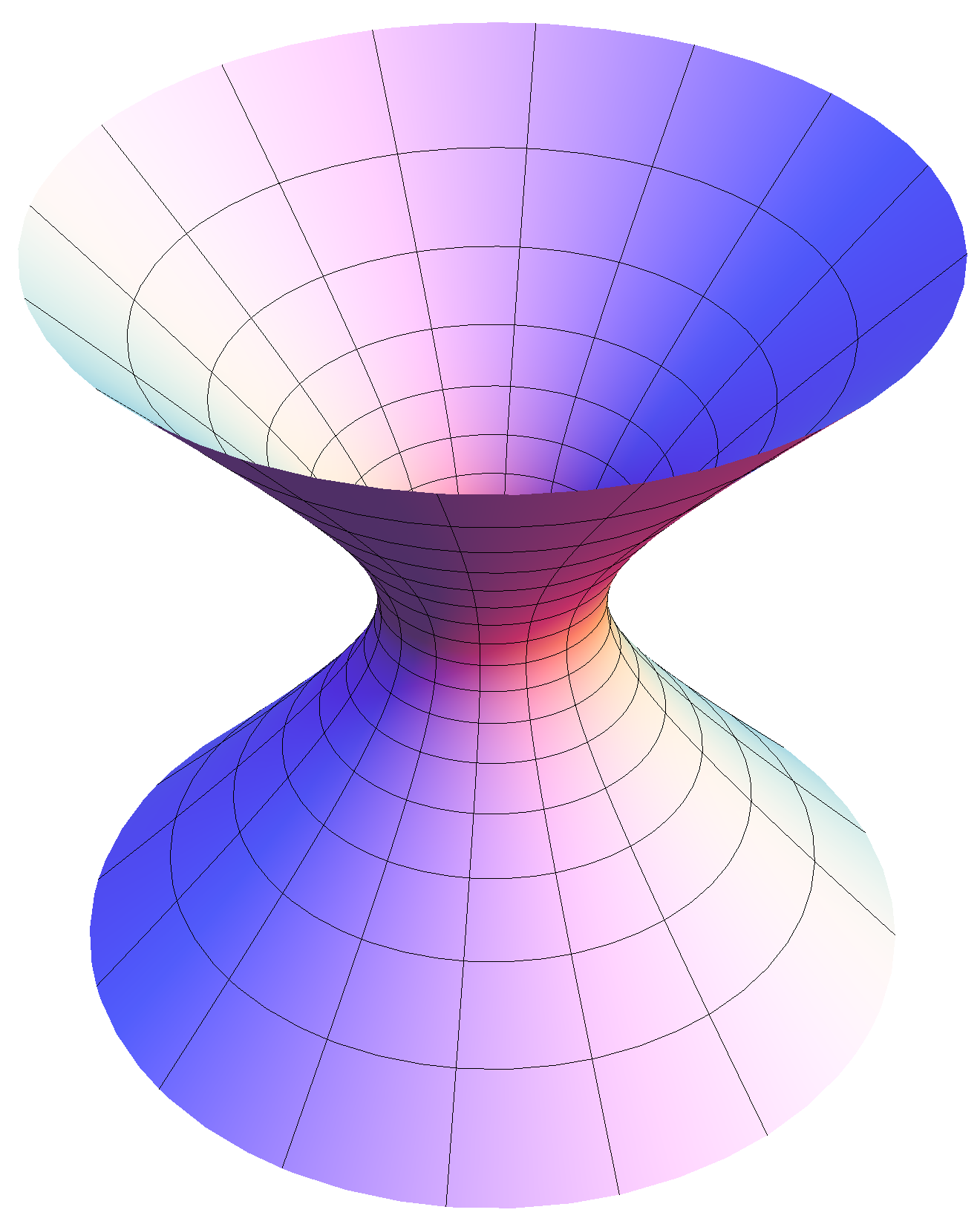} $\qquad$ $\qquad$ $\qquad$
\includegraphics[height=4cm,bb= 0 0 1832 2599]{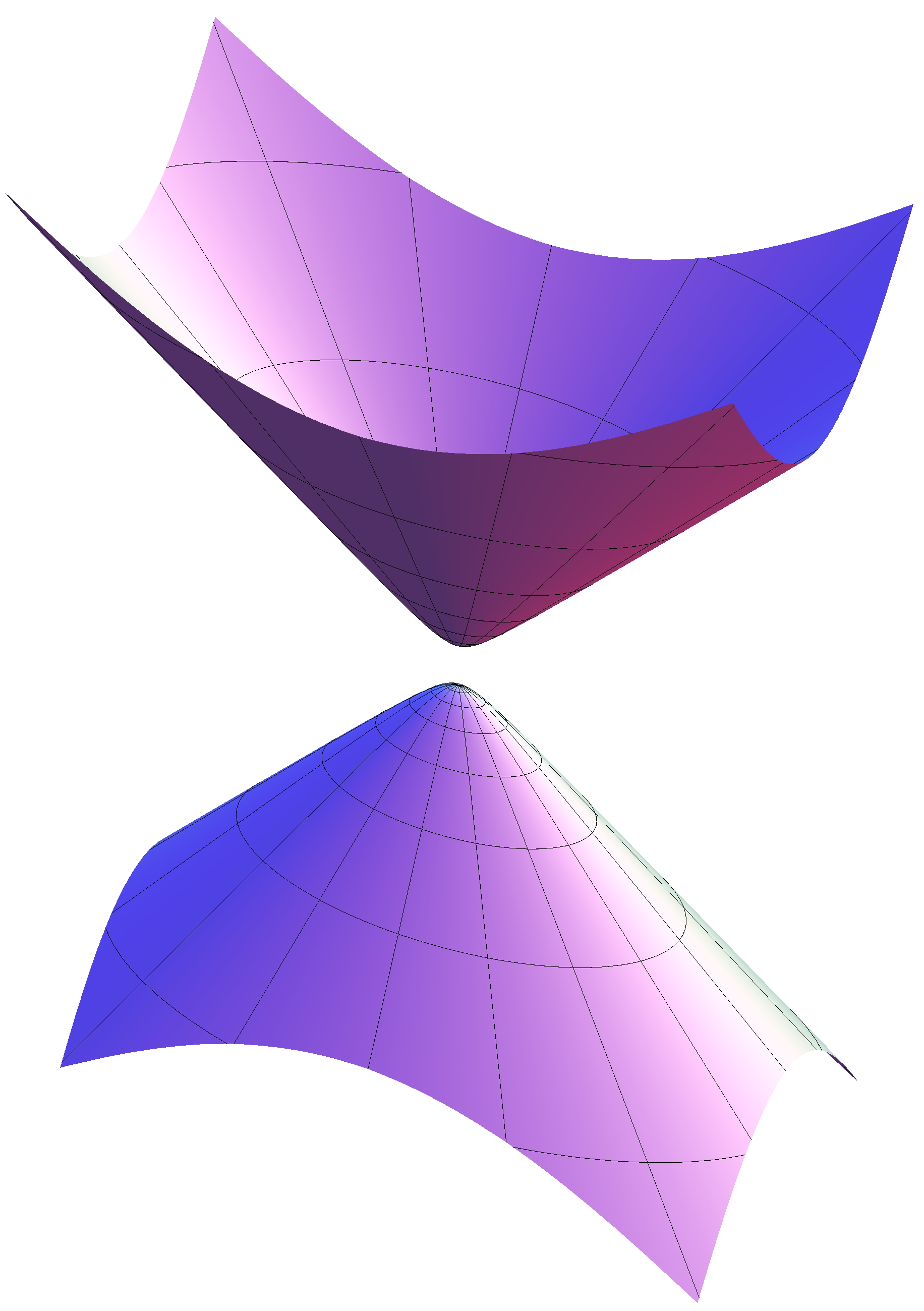}
\caption{Jednokrilni i dvokrilni hiperboloid}\label{fig:hiperboloidi}
\end{figure}

Generalnije, $\Sigma$ defini\v se regularnu povr\v s ako, za sva
$(x,y,z)\in\Sigma$, $$
   \nabla F(x,y,z)\neq0.
$$
\begin{primjer}[Ravan]
   Jedna\v cina ravni kroz tri nekolinearne ta\v cke
   $P_0,P_1,P_2\in\R^3$:
      $$\{P(x,y,z)\ |\ (P-P_0)\cdot (P_1-P_0)\times (P_2-P_0)=0\}.$$
\end{primjer}
\begin{primjer}[Sfera]
   Jedini\v cna sfera, slika \ref{fig:sfera},
      $$\{(x,y,z)\,|\,x^2+y^2+z^2=1\}$$
   je regularna povr\v s: sa $F(x,y,z)=x^2+y^2+z^2-1$ $$
      ||\nabla F(x,y,z)||^2
      = ||2(x,y,z)||^2
      = 4(x^2 + y^2 + z^2) =
       -4 \neq 0.
   $$
\end{primjer}
\begin{primjer}[Jednokrilni hiperboloid]
   Regularna parametrizacija od
      $$\{(x,y,z)\,|\,\left({x\over a}\right)^2+\left({y\over b}\right)^2-\left({z\over c}\right)^2=1\}$$
   je data sa $$
      (u,v)\mapsto(a\cosh u\cos v,b\cosh u\sin v,c\sinh u).
   $$ Vidi sliku \ref{fig:hiperboloidi} lijevo.
\end{primjer}
\begin{primjer}[Dvokrilni hiperboloid]
   Kako bismo pokrili dvostrani hiperboliod, trebamo dvije (regularne) parametrizacije
   jer
      $$\{(x,y,z)\,|\,\left({x\over a}\right)^2+\left({y\over b}\right)^2-\left({z\over c}\right)^2=-1\}$$
   ima dvije konektovane komponente:
   \begin{eqnarray*}
   (u,v)&\mapsto&(a\sinh u\cos v,b\sinh u\sin v,\pm c\cosh u).
   \end{eqnarray*} Vidi sliku \ref{fig:hiperboloidi} desno.
\end{primjer}
\begin{primjer}[Elipsoid]
Iako je elipsoid
\addnumberedpicture{l}{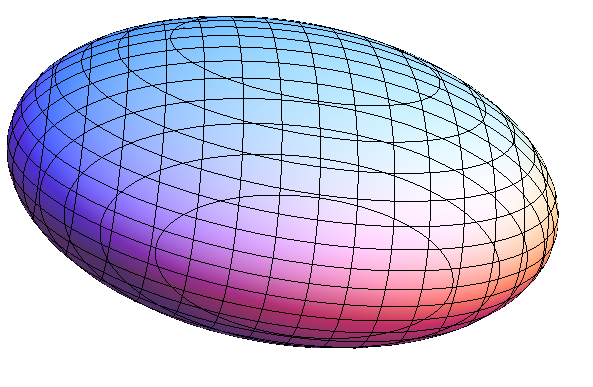}{0.15, bb = 0 0 593 387}{Elipsoid}{fig:elipsoid}
      $$\{(x,y,z)\,|\,\left({x\over a}\right)^2+\left({y\over b}\right)^2+\left({z\over c}\right)^2=1\}$$
   konektovan, nema (kao u slu\v caju sfere) {\it globalne\/} regularne
   parametrizacije (iz istog razloga); trebamo nekoliko regularnih parametrizacija
   kako bismo pokrili elipsoid (vidi vje\v zbu \ref{vjezba:elipsoid}).
\end{primjer}
\begin{primjer}[Elipti\v cna kupa]
Ova povr\v s je parameterizovana sa
\addnumberedpicture{r}{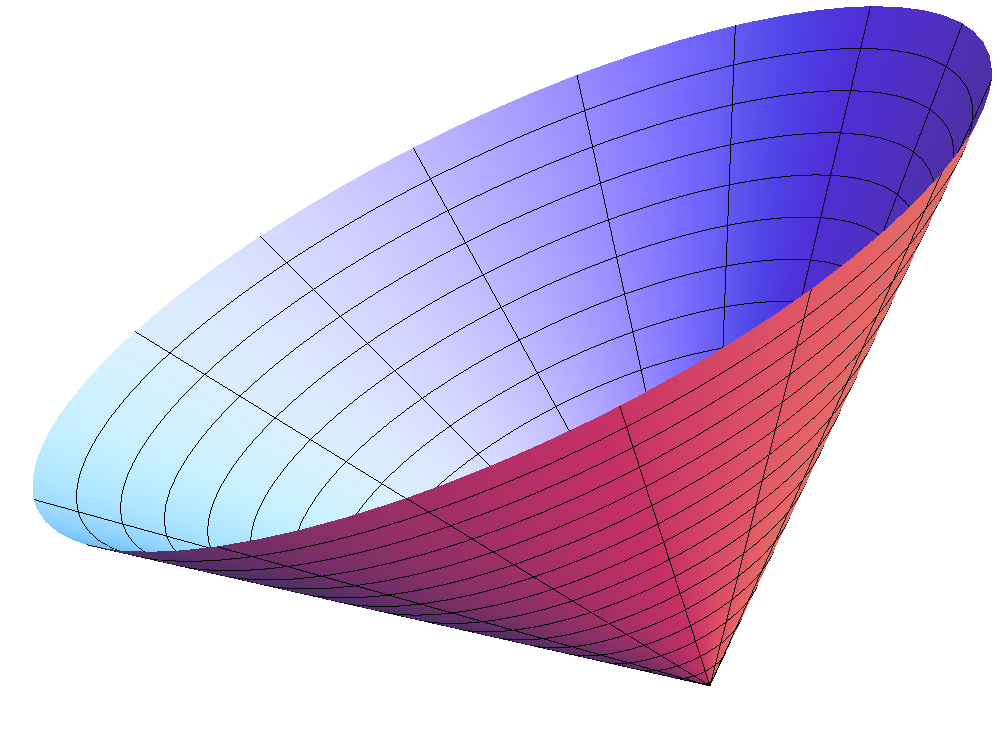}{0.15, bb= 0 0 1000 730}{Elipti\v cna kupa}{fig:kupa}
      $$\{(x,y,z)\,|\,\left({x\over a}\right)^2+\left({y\over b}\right)^2=\left({z\over c}\right)^2\}$$
   i ima singularnost u vrhu kupe; stoga elipti\v cna kupa
   nije regularna povr\v s ali nam slijede\' ca parametrizacija $$
      (u,v)\mapsto(a\,u\cos v,b\,u\sin v,c\,u).
   $$ daje regularnost. Kupa je prastara povr\v s - to je druga povr\v s o kojoj pri\v ca Euklid.
    Euclidovi Elementi, knjiga 11, definicije 18-20 ka\v zu:
   {\it \begin{enumerate}
   \item[18.] Ako jedan krak pravog ugla (jedna kateta) pravouglog trougla ostaje nepokretan, a trougao se oko te prave obr\' ce i vrati u polo\v zaj iz kojeg je poceo kretanje, obuhvacena figura je konus (kupa). Ako je nepokretan krak pravog ugla jednak drugom kraku tog ugla, koji se obrce, konus je pravougli, ako je manji - tupougli, a ako je veci - o\v strougli.
\item[19.] Osa je konusa nepokretna prava oko koje se trougao obrce.
\item[20.] Osnova je konusa krug koji opisuje pokretna du\v z.
   \end{enumerate}}
\end{primjer}
\begin{primjer}[Povr\v s revolucije]
Ako je $t\mapsto \alpha(t)$, $\alpha(t) = (x(t),y(t))$ \\
\addnumberedpicture{r}{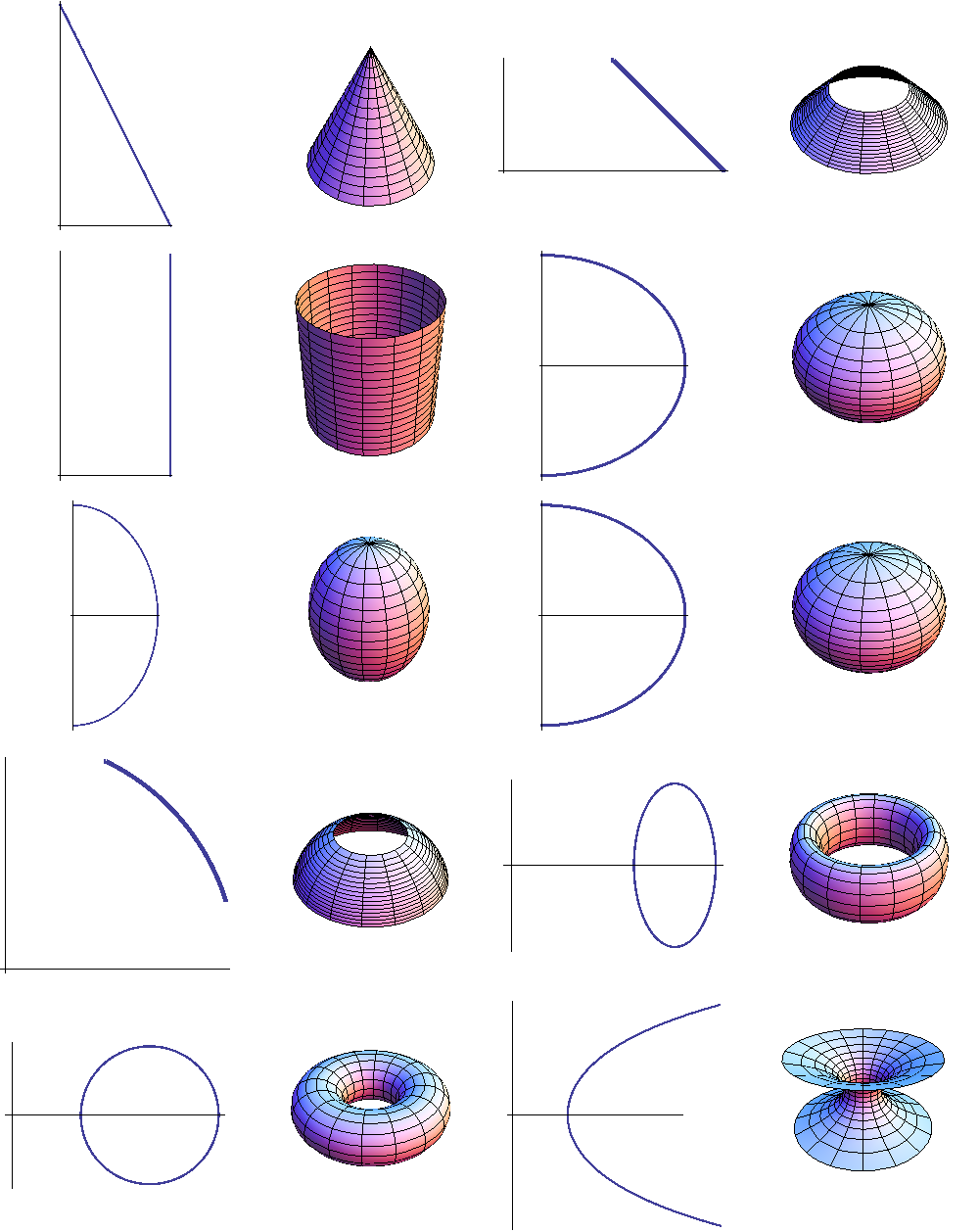}{0.15, bb = 0 0 1000 1258}{Povr\v si revolucije}{fig:revolucija}
regularna kriva u $\mathbb{R}^2$ i neka je bez gubitka op\' cenitosti $x(t)\ne 0$, onda se jedna\v cina \emph{povr\v si revolucije} $\x(u,v)$ dobija rotacijom krive $\alpha(t)$ oko $z$-ose. Njena parametrizacija je data sa
\[
\x(u,v) = (x(u)\cos v, x(u) \sin v, y(u)).
\]
Sa obzirom da je
\[
||\x_u \times \x_v|| = x \sqrt{x'^2+y'^2} \ne 0,
\] imamo da je povr\v s revolucije regularna povr\v s.

Mnoge standardne povr\v si su povr\v si revolucije, vidi sliku \ref{fig:revolucija}.
\end{primjer}
\subsection*{Vje\v zbe}
\begin{vjezba}
Poka\v zite da $\R^2\ni\left(u,v\right)\mapsto\x\left(u,v\right)=\left(u,v,g\left(u,v\right)\right)\in\R^3$, nije \emph{povr\v s}, iako je
$$
   g\left(u,v\right) := \left\{\begin{array}{cl}
      0 & \textrm{ako je } u=v=0  \\
      \frac{uv^2}{u^2+v^2} & \textrm{ako je } \left(u,v\right)\neq\left(0,0\right),
      \end{array}\right.
$$ neprekidna i svi su joj parcijalni izvodi definisani.
\end{vjezba}
\begin{vjezba}
Poka\v zite da $\R^2\ni\left(u,v\right)\mapsto\x\left(u,v\right)=\left(u,v,g\left(u,v\right)\right)\in\R^3$, nije \emph{povr\v s}, iako je
$$
   g\left(u,v\right) := \left\{\begin{array}{cl}
      0 & \textrm{ako je } u=v=0  \\
      \frac{uv^2\sqrt{u^2+v^2}}{u^2+v^4} & \textrm{ako je } \left(u,v\right)\neq\left(0,0\right),
      \end{array}\right.
$$
neprekidna i svi su joj izvodi u pravcu definisani.
\end{vjezba}
\begin{vjezba}
Neka je data implicitna povr\v s $\Sigma:=\{(x,y,z)\,|\,F(x,y,z)=0\}$ i pretpostavimo da je $\Sigma$ regularna povr\v s.
Poka\v zite da je $\displaystyle \n={\nabla F\over|\nabla F|}$ jedini\v cno normalno vektorsko polje na $\Sigma$.
\end{vjezba}
\begin{vjezba}
\begin{enumerate}
\item[(a)] Poka\v zite da je hiperboli\v cni paraboloid
 $$
 \Sigma_1 := \left\{ (x,y,z) \ | \ \frac{z}{c} = \frac{x^2}{a^2} - \frac{y^2}{b^2}\right\}
 $$ regularna povr\v s i na\dj{}ite regularnu parametrizaciju.
\item[(b)]
Poka\v zite da je elipti\v cni paraboloid
$$
\Sigma_2 := \left\{ (x,y,z) \ | \  \frac{z}{c} = \frac{x^2}{a^2} + \frac{y^2}{b^2}\right\}
$$ regularna povr\v s i na\dj{}ite regularnu parametrizaciju.
Zatim ga tako\dj{}er parametrizirajte kao povr\v s revolucije u slu\v caju da je $a=b=c=1$. Zatim tu parametrizaciju iskoristite za nala\v zenje op\' ce parametrizacije.
\end{enumerate}
\end{vjezba}
\begin{vjezba}\label{vjezba:elipsoid}
Poka\v zite da je elipsoid regularna povr\v s i na\dj{}ite regularne parametrizacije koje ga pokrivaju:
\begin{enumerate}
\item[(a)] Koriste\' ci se implicitnim oblikom povr\v si.
\item[(b)] Pode\v savanjem sfernih koordinata na "elipsoidne".
\end{enumerate}
\end{vjezba}
\begin{vjezba}
Posmatrajmo $\Sigma:=\{\left(x,y,z\right)\,|\,\left(\sqrt{x^2+y^2}-R\right)^2+z^2=r^2\}$, gdje je $0<r<R$.
Poka\v zite da je $\Sigma$ regularna povr\v s i na\dj{}ite njenu regularnu parametrizaciju. Kako se zove ova povr\v s?
\end{vjezba}
\begin{vjezba}
\begin{enumerate}
\item[(a)] U Wolfram Mathematica, napravite funkciju \\
\verb|regularna[x_,u_,v_]|
koja ispituje regularnost parametrizovane krive.
\item[(b)] Napravite funkciju
\verb|regularnaImp[F_,x_,y_,z_]|
koja ispituje regularnost implicitno zadate krive.
\item[(c)] Napravite funkciju
\verb|povrsRevolucije[alpha_,u_,v_]|
koja vra\' ca parametrizovanu povr\v s revolucije kreiranu od krive $\alpha(t)$ rotirane oko $z$ ose.
\end{enumerate}
\end{vjezba}

\section{Krive na povr\v sima}
\emph{U punoj verziji teksta.}

\section{Prva fundamentalna forma}
\emph{U punoj verziji teksta.}

\subsection{Normala, tangentna ravan}
\emph{U punoj verziji teksta.}

\subsection{Gaussovo preslikavanje}
\emph{U punoj verziji teksta.}

\subsection*{Vje\v zbe}
\begin{vjezba}
Doka\v zite da je ugao presjeka parametarskih linija na proizvoljnoj povr\v si
\[
    \cos \theta = \frac{F}{\sqrt{EG}}, \qquad
    \sin \theta = \frac{\sqrt{EG-F^2}}{\sqrt{EG}}.
    \]
\end{vjezba}
\begin{vjezba}
Izra\v cunajte prvu fundamentalnu formu, Gaussovo preslikavanje i tangentnu ravan za povr\v s $(u,v)\mapsto \x(u,v)$ datu sa
\[
\x(u,v)=\left(
\frac{\cos v}{\cosh u}, \frac{\sin v}{\cosh u}, \tanh u
\right).
\] Kakva je ovo parametrizacija? Koja je ovo povr\v s?
\end{vjezba}
\begin{vjezba}
Poka\v zite da je kupa $(r,\theta)\mapsto r\,\gamma(\theta)$, gdje je
$\theta\mapsto\gamma(\theta)\in S^2$ du\v zinom luka parametrizovana
sferi\v cna kriva, izometri\v cna ravni.

Poka\v zite da je $\Sigma=\{(x,y,z)\,|\,x^2+y^2=z^2,z>0\}$ izometri\v cna ravni.
\end{vjezba}
\begin{vjezba}
Izra\v cunajte Gaussovo preslikavanje i tangentnu ravan u proizvoljnoj ta\v cki svih do sada predstavljenih parametrizovanih povr\v si.
\end{vjezba}
\begin{vjezba}
Potvrdite da je inverzna  stereografska projekcija sfere tako\dj{}er konformalna parametrizacija.
\end{vjezba}

\begin{vjezba}
\begin{enumerate}
\item[(a)] Napravite  u Wolfram Mathematica animaciju koja pokazuje da Gaussovo preslikavanje mijenja smjer kada jednom obi\dj{}emo M\" obiusovu traku.
\item[(b)] Napravite funkciju \verb|tangentnaRavan[x_,u_,v_,u0_,v0_]| koja ra\v cuna tangentnu ravan povr\v si $\x(u,v)$ u ta\v cki $\x(u0,v0)$.
\item[(c)] Napravite funkciju \verb|gaussovo[x_,u_,v_]| koja ra\v cuna Gaussovo preslikavanje date povr\v si.
\item[(d)] Napravite funkciju \verb|pff[x_,u_,v_]| koja ra\v cuna prvu fundamentalnu formu date povr\v si i vra\' ca je kao listu \verb|{E,F,G}|.
\end{enumerate}
\end{vjezba}

\section{Linijske i razvojne povr\v si}
\emph{U punoj verziji teksta.}

\subsection*{Vje\v zbe}
\begin{vjezba}
Poka\v zite da je jednokrilni hiperboloid $\Sigma=\{(x,y,z)\,|\,x^2+y^2=1+z^2\}$
linijska povr\v s te ispitajte kako se njegovo Gaussovo preslikavanje mijenja du\v z pravih linija na $\Sigma$.
\end{vjezba}

\begin{vjezba}
Neka je $s\mapsto\gamma(s)$ du\v zinom luka parametrizovana kriva sa
$\kappa\neq0$.
\begin{enumerate}
\item[(a)] Izra\v cunajte prvu fundamentalnu formu odgovaraju\' ce tangentne razvojne povr\v si.
\item[(b)] Na\dj{}ite ortogonalnu ($F\equiv0$) reparametrizaciju povr\v si (kao razvojne povr\v si).
\end{enumerate}
\end{vjezba}

\begin{vjezba}
Neka je $(u,v)\mapsto\x(u,v)=\gamma(u)+v\,\eta(u)$ razvojna povr\v s.
Poka\v zite da je $v\mapsto\n(u,v)$ paralelno normalno polje du\v z
$v\mapsto\x(u,v)$.
\end{vjezba}

\begin{vjezba}
Neka je $(u,v)\mapsto\x(u,v)=\gamma(u)+v\,\eta(u)$ linijska povr\v s.
Poka\v zite da su linije $v\mapsto\x(u,v)$ ($u$ fiksno), asimptotske linije.
\end{vjezba}

\section{Druga fundamentalna forma}
\emph{U punoj verziji teksta.}

\subsection*{Vje\v zbe}
\begin{vjezba}
Izra\v cunajte druge fundamentalne forme sfere, helikoida i jednokrilnog hiperboloida.
\end{vjezba}

\begin{vjezba}
\begin{enumerate}
\item[(a)] Napravite  u Wolfram Mathematica animacije koje pokazuju da su kupa, hiperboloid, kru\v zni helikoid linijske povr\v si.
\item[(b)] Napravite funkciju \verb|dff[x_,u_,v_]| koja ra\v cuna prvu fundamentalnu formu date povr\v si i vra\' ca je kao listu \verb|{e,f,g}|.
\end{enumerate}
\end{vjezba}

\subsection{Meusnierova teorma}
\emph{U punoj verziji teksta.}

\begin{teorem}[Meusnier]
Oskulatorne kru\v znice presje\v cnih krivih povr\v si $\Sigma$ i ravni koja
sadr\v zi datu ne-asimptotski tangentni pravac povr\v si $\Sigma$
u $(x,y,z)\in\Sigma$ su sadr\v zani na sferi koja doti\v ce $\Sigma$ u $(x,y,z)$.
\end{teorem}
\begin{dokaz}
Fiksirajmo tangentni pravac $T$ u ta\v cki $(x,y,z)\in\Sigma$,
tojest $T\perp\n(x,y,z)$, gdje $\n(x,y,z)$ ozna\v cava jedini\v cni normalni vektor
u $(x,y,z)$.  Neka je $\kappa_n$ normalna krivina pravca
$T$; jer za $T$ pretpostavljamo da je ne--asimptotski pravac
imamo $\kappa_n\neq0$ i mo\v zemo posmatrati ``radijus krivine''
$R_n={1\over\kappa_n}$ (bez gubitka generalnosti mo\v zemo pretpostaviti da je $R_n>0$ nakon
\v sto mogu\' ce zamjenimo $\n(x,y,z)$ sa $-\n(x,y,z)$).
Neka je $c_n:=(x,y,z)+R_n\,\n(x,y,z)$.

Sada, fiksirajmo jedini\v cni vektor $N\perp T$, koji {\it nije\/} tangentni
vektor ($N\not\perp\n(x,y,z)$), i posmatrajmo krivu presjeka ravni
$$
   {\cal E} := \{(x,y,z)+\lambda T+\mu N\,|\,\lambda,\mu\in\R\}
$$
sa $\Sigma$.
Pretpostavi\' cemo da je $N$ izabrana tako da $\cos\varphi:=N\cdot\n(x,y,z)>0$.
Jer je $\kappa={1\over\cos\varphi}\kappa_n$ dobivamo radijus i centar oskulatorne
kru\v znice presje\v cne krive:$$
   R = \cos\varphi\,R_n
      \quad{\textrm{ i }}\quad
   c = (x,y,z) + R\,N.
$$
Stoga oskulatorna kru\v znica mo\v ze biti parametrizovana sa$$
   t \mapsto \alpha(t) = c + R\,\{\cos t\,T + \sin t\,N\};
$$
onda
\begin{eqnarray*}
   |\alpha(t) - c_n|^2
   &=& |R\cos t\,T+R(1+\sin t)\,N+R_n\n(x,y,z)|^2 \\
   &=& 2R(R-R_n\cos\varphi)(1+\sin t)+R_n^2 \\
   &=& R_n^2
\end{eqnarray*}
\v sto pokazuje da je oskulatorna kru\v znica na sferi sa centrom $c_n$
i radijusom $R_n$; primjetite da ova sfera doti\v ce $\Sigma$ u $(x,y,z)$
jer je $(x,y,z)$ na njoj i ima istu normalu $\n(x,y,z)$ kao
povr\v s u ovoj ta\v cki, tako da se tangentne ravni povr\v si i sfere
podudaraju.
\end{dokaz}

\section{Weingartenov tenzor}
\emph{U punoj verziji teksta.}

\subsection{Krivine na povr\v si}
\emph{U punoj verziji teksta.}

\subsection*{Vje\v zbe}
\begin{vjezba}
Izra\v cunajte Gaussovu krivinu razvojne povr\v si.
\end{vjezba}
\begin{vjezba}
Izra\v cunajte Gaussovu, srednju i principalne krivine  sfere, helikoida i jednokrilnog hiperboloida.
\end{vjezba}

\begin{vjezba}
Neka je $(u,v)\mapsto \sigma(u,v)$ regularna parametrizovana povr\v s.
\begin{itemize}
\item[(a)] Doka\v zite da su $\n(u,v)$ i $\n_u\times\n_v (u,v)$ me\dj{}usobno paralelni vektori za sva $u,v$.
\item[(b)]  Izra\v cunajte $(\n_u\times \n_v).\n$ i zaklju\v cite da je
$$
\n_u\times \n_v = K\sigma_u \times \sigma_v
$$
\item[(c)] Neka nam je sada data nova povr\v s: $(u,v)\mapsto \rho(u,v) + a \n(u,v)$, $a\in \R$.
Doka\v zite da su $\n(u,v)$ i $\rho_u\times\rho_v (u,v)$ me\dj{}usobno paralelni vektori za sva $u,v$. Doka\v zite da je
$$
\rho_u\times\rho_v = (1-2H\ a + K\ a^2) \sigma_u \times \sigma_v
$$
\end{itemize}
\end{vjezba}

\begin{vjezba}
Pretpostavimo da povr\v s bez umbili\v cnih ta\v caka $(u,v)\mapsto\x(u,v)$ ima konstantnu
srednju krivinu $H\neq0$.

Doka\v zite da povr\v s $(u,v)\mapsto\x^\ast(u,v):=\x(u,v)+{1\over H}\n(u,v)$ ima
konformalno ekvivalentnu metriku $\fff^\ast={H^2-K\over H^2}\,\fff$
i konstantnu srednju krivinu $H^\ast=H$.
\end{vjezba}

\section{Eulerov teorem}
\emph{U punoj verziji teksta.}

\subsection*{Vje\v zbe}
\begin{vjezba}
Doka\v zite Eulerov teorem, odnosno `popunite praznine' dokaza iznad.
\end{vjezba}
\begin{vjezba}
\begin{enumerate}
\item Koriste\' ci se ve\' c izra\v cunatim fundamentalnim formama sa prethodnih vje\v zbi, odredite povr\v si i ta\v cke na tim povr\v sima koje su elipti\v cne, hiperboli\v cne i paraboli\v cne.
\item[(c)] Demonstrirajte kako se priroda tih ta\v caka ogleda u presjeku tangentnih ravni sa povr\v si, koriste\' ci se Wolfram Mathematicom.
\end{enumerate}
\end{vjezba}

\section{Minimalne povr\v si}
\emph{U punoj verziji teksta.}

\subsection*{Vje\v zbe}
\begin{vjezba}
Poka\v zite da je helikoid minimalna povr\v s i odredite njegove asimptotske linije
i linije krivine.
\end{vjezba}
\begin{vjezba}
Poka\v zite da je {\it Enneperova povr\v s\/}
$(u,v)\mapsto(u^3-3u(1+v^2),v^3-3v(1+u^2),3(u^2-v^2))$
konformalno parametrizovana minimalna povr\v s.
\end{vjezba}
\begin{vjezba}
Neka je $(u,v)\mapsto\x(u,v)$ konformalno parametrizovana minimalna povr\v s.
Poka\v zite da se njena druga fundamentalna forma mo\v ze napisati kao
$\sff=\textrm{Re }\{(e-if)(du+idv)^2\}$ i da je
$(e^\alpha-if^\alpha)=e^{i\alpha}(e-if)$ za
povr\v si $\x^\alpha$ asocijativne porodice od $\x$.

Zaklju\v cite da asimptotske linije i linije krivine povr\v si $\x$ postaju
asimptotske linije i linije krivine povr\v si $\x^\ast$, respektivno.
\end{vjezba}
\begin{vjezba}
Poka\v zite da je povr\v s revolucije krive $t\mapsto \cosh t$ (katenoida) minimalna.
\end{vjezba}

~ \newpage ~ \hspace{1cm} ~  \par ~ \\

\vfill {\scalebox{0.7}{{\stechak lkkkkkkkkkkkkkkkkkd}}}

\chapter{Fundamentalne jedna\v cine povr\v si}
\pagestyle{main2}
\minitoc
\vspace{1cm}
Ovo poglavlje ima dvostruku svrhu. Prvo, uvest \' cemo pojam ``kovarijantnog izvoda'' ili
``Levi-Civita konekcije'' i vidjet \' cemo da ona samo zavisi o prvoj
fundamentalnoj formi povr\v si;
drugo, postavit \' cemo scenu za formulaciju Gauss-Codazzijevih jedna\v cina
u jednoj od slijede\' cih sekcija, tako \v sto \' cemo zapisati strukturne
jedna\v cine za (ne-ortonormirani) okvir povr\v si.
\section{Kovarijantni izvod i jedna\v cine Gauss -- Weingartena}
U ovom poglavlju radi kratko\' ce koristit \' cemo ``tenzorsku notaciju'', koja se na prvi pogled mo\v ze \v ciniti komplikacijom. Me\dj{}utim, vrlo brzo \' cemo se uvjeriti u njezine prednosti.
\begin{notacija}
Parcijalne izvode ozna\v cavat \' cemo sa
$$
   \partial_1 = \frac{\partial}{\partial u}
      \quad{\textrm{ i }}\quad
   \partial_2 = \frac{\partial}{\partial v},
$$
a fundamentalne forme sa
$$
   \fff = \left(\begin{array}{cc}g_{11} & g_{12} \\ g_{21} & g_{22}\end{array}\right)
      \quad{\textrm{ i }}\quad
   \sff = \left(\begin{array}{cc}h_{11} & h_{12} \\ h_{21} & h_{22}\end{array}\right).
$$
\end{notacija}
Osnovna ideja je da formuli\v semo strukturne jedna\v cine za
okvir $F=(\x_u,\x_v,\n)$: kako sada imamo dvodimenzionalnu domenu
ovi imaju slijede\' ci oblik  $$
   F_u = F\,\Phi \quad{\textrm{ i }}\quad F_v = F\,\Psi
$$
sa dvije funkcije sa matri\v cnim vrijednostima $\Phi$ i $\Psi$.
Mi ne\' cemo, me\dj{}utim, koristiti ovaj formalizam.

\begin{defn}[Kovarijantni izvod]
Neka je $(u,v)\mapsto\xi(u,v)$ tangencijalno vektorsko polje du\v z parametrizovane
povr\v si $(u,v)\mapsto\x(u,v)$;  \emph{kovarijantni izvod} tog polja
u $(u,v)$ u pravcu $(\lambda,\mu)$ je
$$
   (\nabla_{(\lambda,\mu)}\xi)|_{(u,v)}
   :=  \{(\lambda\xi_u+\mu\xi_v)  - ((\lambda\xi_u+\mu\xi_v)  \cdot\n)\, \n\} (u,v).
$$
$\nabla$ se tako\dj{}er zove \emph{Levi-Civita konekcija}.

Koeficijenti $\Gamma_{ij}^k$ u
$$
   \nabla_j(\partial_i\x) :
   =   \partial_j\partial_i\x -   (\partial_j\partial_i\x\cdot\n)\n
   = \sum_k\Gamma_{ij}^k\,\partial_k\x
$$
se zovu \emph{Christoffelovi simboli}.
\end{defn}
Imamo da je
\begin{eqnarray*}
   \nabla_{(1,0)}\x_u =: \Gamma_{11}^2\x_u + \Gamma_{11}^2\x_v,  &
   \nabla_{(0,1)}\x_u =: \Gamma_{12}^2\x_u + \Gamma_{12}^2\x_v; \\
   \nabla_{(1,0)}\x_v =: \Gamma_{21}^2\x_u + \Gamma_{21}^2\x_v,  &
   \nabla_{(0,1)}\x_v =: \Gamma_{22}^2\x_u + \Gamma_{22}^2\x_v.
\end{eqnarray*}
\begin{lema} Christoffelovi simboli su simetri\v cni, odnosno
$$\Gamma_{ij}^k=\Gamma_{ji}^k.$$
\end{lema}
\begin{dokaz}
$\partial_2\partial_1\x=\x_{uv}=\x_{vu}=\partial_1\partial_2\x$.
\end{dokaz}
\begin{remark}
Sada mo\v zemo formulisati strukturne jedna\v cine za (ne-ortonormirani)
adaptirani okvir $(u,v)\mapsto F(u,v):=(\x_u,\x_v,\n)(u,v)$:
$$
   \partial_j\partial_i\x = \sum_k\Gamma_{ij}^k\,\partial_k\x + h_{ij}\n
$$
ili ekvivalentno,
\begin{eqnarray}\label{struk1}
   \x_{uu} &=& \Gamma_{11}^1\x_u + \Gamma_{11}^2\x_v + e\n\,\\
   \label{struk2}
   \x_{uv} &=& \Gamma_{12}^1\x_u + \Gamma_{12}^2\x_v + f\n\,\\
   \label{struk3}
   \x_{vv} &=& \Gamma_{22}^1\x_u + \Gamma_{22}^2\x_v + g\n\,
\end{eqnarray}
i drugo (\v sto dobijamo iz $d_{(u,v)}\n=-d_{(u,v)}\x\circ\wtf|_{(u,v)}$)
$$
   \partial_i\n = -\sum_{k,l}s_{ki}\partial_k\x,
$$
gdje su $s_{ki}=\sum_lg^{kl}h_{li}$  koeficijenti $\wtf$, ili ekvivalentno,
\begin{eqnarray*}
   \n_u &=& -{1\over EG-F^2}\{(Ge-Ff)\,\x_u+(Ef-Fe)\,\x_v\} \\
   \n_v &=& -{1\over EG-F^2}\{(Gf-Fg)\,\x_u+(Eg-Ff)\,\x_v\}
\end{eqnarray*}
Ovo skupa \v cini ``{\it jedna\v cine Gauss-Weingartena\/}''.
\end{remark}
\begin{primjer}
\v Zelimo izra\v cunati eksplicitne formule za Christoffelove simbole $\Gamma^{k}_{ij}$ koriste\' ci se koeficijentima prve fundamentalne forme.

Koriste\' ci se strukturnim jedna\v cinama (\ref{struk1}),(\ref{struk2}),(\ref{struk3}), te
uzimaju\' ci skalarni proizvod svake jedna\v cine sa  $x_u$ i $x_v$, dobivamo
\begin{align*}
\Gamma^{1}_{11} E + \Gamma^{2}_{11} F &= \frac12 E_u, &  \Gamma^{1}_{11} F + \Gamma^{2}_{11} G &= F_u-\frac12 E_v, \\
\Gamma^{1}_{12} E + \Gamma^{2}_{12} F &= \frac12 E_u, &  \Gamma^{1}_{12} F + \Gamma^{2}_{12} G &= \frac12 G_u, \\
\Gamma^{1}_{22} E + \Gamma^{2}_{22} F &= F_v - \frac12 G_u, &  \Gamma^{1}_{22} F + \Gamma^{2}_{22} G &= \frac12 G_v.
\end{align*}
Rje\v savaju\' ci ovaj sistem u parovima, dobijamo
\begin{align*}
\Gamma^{1}_{11} &= \frac{GE_u-2FF_u+FE_v}{2(EG-F^2)}, & \Gamma^{2}_{11} &= \frac{2EF_u-EE_v-FE_u}{2(EG-F^2)}, \\
\Gamma^{1}_{12} &= \frac{GE_v-FG_u}{2(EG-F^2)}, &
\Gamma^{2}_{12} &= \frac{EG_u-FE_v}{2(EG-F^2)}, \\
\Gamma^{1}_{22} &= \frac{2GF_v-GG_u-FG_v}{2(EG-F^2)}, &
\Gamma^{2}_{22} &= \frac{EG_v - 2FF_v + FG_u}{2(EG-F^2)}.
\end{align*}
\end{primjer}
\begin{lema}[Koszul]
$$
   \sum_mg_{km}\Gamma_{ij}^m
   = {1\over2}\{\partial_ig_{jk}+\partial_jg_{ik}-\partial_kg_{ij}\}.
$$
\end{lema}
\begin{dokaz}
Ra\v cunamo s desna na lijevo:
$$\textstyle
   \partial_ig_{jk}+\partial_jg_{ik}-\partial_kg_{ij}
   = 2\partial_k\x\cdot\partial_j\partial_i\x
   = 2\sum_mg_{km}\Gamma_{ij}^m,
$$ \v sto nam daje rezultat.
\end{dokaz}

Kako je $\fff$ pozitivno definitna, Koszulov identitet mo\v ze biti rije\v sen za
$\Gamma_{ij}^m$ i dobijamo:
\begin{cor}
Kovarijantni izvod $\nabla$ samo zavisi o $\fff$.
\end{cor}
\begin{cor}
Christofellovi simboli mogu se napisati kao
$$
\Gamma^m_{ij} = \frac12 \sum_k g^{km}(\partial_j g_{ik} + \partial_i g_{jk}-\partial_k g_{ij}).
$$
\end{cor}

\section{Gaussova Theorema egregium i jedna\v cine Codazzija}
\emph{U punoj verziji teksta.}

\subsection{Primjene jedna\v cina Gaussa i Codazzija}
\emph{U punoj verziji teksta.}

\subsection{Totalno umbili\v cne povr\v si}
\emph{U punoj verziji teksta.}

\section{Fundamentalna teorema teorije povr\v si}
\emph{U punoj verziji teksta.}

\subsection*{Vje\v zbe}
\begin{vjezba}
 Dokazati posljedicu s predavanja :
$$
\Gamma^m_{ij} = \frac12 \sum_k g^{km}(\partial_j g_{ik} + \partial_i g_{jk}-\partial_k g_{ij}).
$$
\end{vjezba}
\begin{vjezba}
Izra\v cunati Christofellove simbole za sferu, helikoid i jednokrilni hiperboloid.
\end{vjezba}
\begin{vjezba}
Napraviti funkciju \verb|Christoffel[x_,u_,v_]| u Wolfram Mathematici.
\end{vjezba}
\begin{vjezba}
Neka je $(u,v)\mapsto\x(u,v)$ paramterizacija linijom krivine
povr\v si sa konstantnom Gaussovom krivinom $K\equiv-1$ tako da, bez gubitka op\' cenitosti,
$\kappa_1=\tan\varphi$ i $\kappa_2=-\cot\varphi$ sa odgovaraju\' com funkcijom $\varphi$.

Poka\v zite da postoji reparametrizacija linijom krivine, $u=u(\tilde u)$
i $v=v(\tilde v)$, tako da $$\textstyle
   \tilde\fff = \cos^2\varphi\,d\tilde u^2 + \sin^2\varphi\,d\tilde v^2
      \quad{\text{ i }}\quad
   \tilde\sff = {1\over2}\sin2\varphi\,(du^2-dv^2).
$$
(Pomo\' c: poka\v zite da je $({E\over\cos^2\varphi})_v=({G\over\sin^2\varphi})_u=0$
iz Codazzijevih jedna\v cina.)
\end{vjezba}
\begin{vjezba}
Neka je $(u,v)\mapsto\x(u,v)$ parametrizacija linijom krivine minimalne povr\v si.
Poka\v zite da $E^2K\equiv const$ i zaklju\v cite da
$\ln E$ mora zadovoljavati Liouvilleovu jedna\v cinu
$$\Delta\ln E=const\,e^{-\ln E}.$$
\end{vjezba}
\begin{vjezba}
Poka\v zite da je Gaussova krivina povr\v si invarijantna pod
reparamterizacijom i zaklju\v cite da Gaussova krivina nestaje
ako je povr\v s (lokalno) izometri\v cna ravni, tojest prima
izometri\v cnu parametrizaciju.
\end{vjezba}

\vfill {\scalebox{0.7}{{\stechak lkkkkkkkkkkkkkkkkkd}}}

\chapter{Geometrija na povr\v si}
\pagestyle{main2}
\minitoc
\vspace{1cm}
Prisjetimo se da je geodezijska krivina $\kappa_g$ du\v zinom luka
parametrizovane krive $s\mapsto\gamma(s)=\x(u(s),v(s))$ definisana sa
$$
   \kappa_g := T'\cdot(N\times T) = |N,T,T'|,
$$
gdje je $N(s)=\n(u(s),v(s))$ povr\v sinska normala koju posmatramo kao
jedini\v cno normalno vektorsko polje du\v z krive $\gamma$.

Kao kod normalne krivine tra\v zimo analiti\v cku reprezentaciju:
\begin{eqnarray*}
\kappa_g
   &=& |N,\gamma',\gamma''| = |\n, \x_uu'+\x_vv', \x_uu''+\x_vv'' +\x_{uu}u'^2+2\x_{uv}u'v'+\x_{vv}v'^2|\\
   &=& |\n,\x_uu'+\x_vv',\x_uu''+\x_vv''| \\
   &+& |\n,\x_uu'+\x_vv',\x_{uu}u'^2+2\x_{uv}u'v'+\x_{vv}v'^2| \\
   &=& \sqrt{EG-F^2}\,
   \left\{
      \left|{u'\atop v'}\,{u''\atop v''}\right|
    + \left|{u'\atop v'}\,{
      (\Gamma_{11}^1u'^2+2\Gamma_{12}^1u'v'+\Gamma_{22}^1v'^2)\atop
      (\Gamma_{11}^2u'^2+2\Gamma_{12}^2u'v'+\Gamma_{22}^2v'^2)
      }\right|\right\}.
\end{eqnarray*}
Iako je ovo vrlo nezgrapan izraz, nau\v cili smo korisnu informaciju, na osnovu rezultata iz prethodnog poglavlja :

\begin{teorem}
Geodezijska krivina $\kappa_g$ samo zavisi o $\fff$.
\end{teorem}

\begin{remark}
Primjetite da je, sli\v cno kao u slu\v caju prostorne krive, geodezijska krivina
$$
   \kappa_g = {|N,\gamma',\gamma''|\over||\gamma'||^3}
$$
u odnosu na {\it bilo koju\/} parametrizaciju;
Sjetite se da je $||\gamma'||^2=Eu'^2+2Fu'v'+Gv'^2$.
\end{remark}

\section{Geodezije}
\emph{U punoj verziji teksta.}

\subsection*{Vje\v zbe}
\begin{vjezba}
Neka je povr\v s $(u,v)\mapsto\x(u,v)$ parametrizirana sa $\x(u,v)=(u \cos v, u\sin v, v)$. Poka\v zite da je kriva
$\beta(t)=\x(t,0)$ geodezija povr\v si $\x$.
\end{vjezba}
\smallskip
\begin{vjezba}
Na\dj{}ite geodezije na ravni ${\cal E}\subset\R^3$.
\end{vjezba}
\smallskip
\begin{vjezba}
Neka je $\Sigma\subset\R^3$ cilindar na krivoj u $(x,y)$-ravni.
Doka\v zite da je geodezija na $\Sigma$ (op\v sti) heliks.
\end{vjezba}
\smallskip
\begin{vjezba}
Parametri\v site $S^2(R)=\{(x,y,z)\,|\,x^2+y^2+z^2=R^2\}$ pomo\' cu geodezijskih polarnih
koordinata oko $P=(0,0,R)$.
\end{vjezba}
\smallskip
\begin{vjezba}
U Mathematici napraviti program koji provjerava da li je kriva geodezija ili ne. Da li mo\v zete napraviti program koji nalazi geodezije?
\end{vjezba}

\vfill {\scalebox{0.7}{{\stechak lkkkkkkkkkkkkkkkkkd}}}

\appendix

\chapter{Diferencijalni ra\v cun realnih funkcija - predznanje}\label{appendix:predznanje}
\pagestyle{appendix2}

\begin{defn}
Preslikavanje $f:\R^m ~^\circ\supset U\to\R^n$ je {\it diferencijabilno} u ta\v cki $p\in U$ ukoliko postoji linearno preslikavanje $d_pf\in{\textrm Hom}(\R^m,\R^n)$ tako da $$
   \lim_{h\rightarrow 0}{f(p+h)-f(p)-d_pf(h)\over|h|} = 0.
$$
\end{defn}
\begin{defn}
Ako je $f$ diferencijabilna u $p\in U\subset^\circ\R^m$, onda je njen izvod
$$
   d_pf \simeq Jf(p) = \left(\begin{array}{ccc}
      {\partial f_1\over\partial x_1}&\cdots&{\partial f_1\over\partial x_m}\\
   \vdots && \vdots \\
      {\partial f_n\over\partial x_1}&\dots&{\partial f_n\over\partial x_m}
      \end{array}\right).
$$
\end{defn}
Primjetite da $f$ ne mora biti diferencijabilna u $p$ ako $Jf(p)$ postoji:
posmatrajmo $$
   f:\R^2\mapsto\R, \
   (u,v)\mapsto f(u,v):=
   \left\{
        \begin{array}{rcl}
            \frac{uv^2\sqrt{u^2+v^2}}{u^2+v^4} & for & (u,v)\neq(0,0) \\
            0 & for & (u,v)=(0,0)
        \end{array}
   \right. .
$$
Ako me\dj{}utim $Jf$ postoji svuda i neprekidna je, onda je $f$
{\it neprekidno diferencijabilna\/} i, konkretno, je diferencijabilna u svakoj ta\v cki.
Ako je $f$ diferencijabilna u $p\in U\subset^\circ\R^m$ onda je njen izvod u pravcu $w\in\R^m$ dat sa $$
   d_pf(w) = Jf(p)\cdot w = \lim_{t\rightarrow 0}\frac{f(p+tw)-f(p)}{t}.
$$
Primjetite da postoje funkcije za koje svi izvodi u pravcu postoje
u ta\v cki, ali koje nisu diferencijabilne u toj ta\v cki.

{\bf Dogovor.}
Za potrebe ovoga predmeta mi \' cemo pretpostavti da je svaka funkcija
diferencijabilna onoliko puta koliko \v zelimo, tj, da je svaka funkcija u $C^\infty$.
\begin{teorem}[Teorem inverznog preslikavanja]
Neka je $f:\R^n ~^\circ\supset U\to\R^n$ neprekidno diferencijabilna i
pretpostavimo da je $d_pf:\R^n\to\R^n$ invertibilna matrica u $p\in U$.
Onda je $f$ lokalno invertibilna funkcija u okolini $p$ sa neprekidno diferencijabilnim inversom; konkretno, $d_{f(p)}(f^{-1})=(d_pf)^{-1}$.
\end{teorem}
\begin{teorem}[Teorem implicitnog preslikavanja] \label{TeoremImplicitnogPreslikavanja}
Neka je $F:\R^m\times\R^k ~^\circ\supset U\times V\to\R^k$ neprekidno
diferencijabilna i pretpostavimo da je $\R^k$-dio od $d_{(p,q)}F$
invertibilna za $(p,q)\in U\times V$.

Onda skup $F(x,y)\equiv F(p,q)$ mo\v zemo napisati kao graf
$y=g(x)$ neprekidno diferencijabilne funkcije $g$.
\end{teorem}
\begin{defn}\label{defn:imerzija}
Neka je $f:\R^m ~^\circ\supset U\to\R^n$ neprekidno diferencijabilna.
Onda se $f$ zove:
\begin{itemize}
\item {\it imerzija\/} ako je $d_pf$ injekcija za sve $p\in U$
   (konkretno, $m\leq n$);
\item {\it submerzija\/} ako je $d_pf$ surjekcija za sve $p\in U$
   (konkretno, $m\geq n$);
\item {\it difeomorfizam\/} ako je $d_pf$ bijekcija
   za sve $p\in U$ (konkretno, $m=n$).
   \end{itemize}
\end{defn}
Obi\v cna diferencijalna jedna\v cina (reda $n$)
je jedna\v cina oblika
$$
   y^{(n)}(x) = f(x,y(x),y'(x),\dots,y^{(n-1)}(x))
\eqno{(\dagger)}$$
za nepoznatu funkciju $y=y(x)$ koja zavisi od realne promjenljive $x$
-- ali koja mo\v ze imati vrijednosti u $\R^m$.

Svaka takva ODJ mo\v ze biti druga\v cije napisana kao (sistem) ODJ reda $n=1$
tako \v sto uvedemo izvode kao nove funkcije:
sa $y_k:=y^{(k-1)}$ jedna\v cina $(\dagger)$ je ekvivalentna sistemu
\begin{eqnarray*}
   y_1'(x) \hfill&=& y_2(x) \hfill\\
   &\cdots&\\
   y_{n-1}'(x) \hfill&=& y_n(x) \hfill\\
   y_n'(x) \hfill&=& f(x,y_1(x),\dots,y_n(x)) \hfill
\end{eqnarray*}
pa stoga nikad ne moramo misliti o ODJ vi\v seg reda.
\begin{teorem}[Picard-Lindel\"of]\label{thm:picard}
Neka je $\R\times\R^n ~^\circ\supset I\times U\ni(x,y)\mapsto f(x,y)\in\R^n$
neprekidna i Lipschitz neprekidna po $y$ i $(x_0,y_0)\in I\times U$;
onda {\it problem po\v cetne vrijednosti\/}
$$
   y'(x) = f(x,y(x)), \quad y(x_0)=y_0
\eqno{(\star)}$$
ima jedinstveno {\footnote{Naravno, ograni\v citi
``maksimalno'' rje\v senje na manji interval daje ``novo'' rje\v senje.
Uporedite ovo sa {\it Peanovom teoremom\/} koja samo zahtjeva neprekidnost
ali ne daje jedinstvenost (obja\v snjenje: $y'=\sqrt{|y|}$, $y(0)=0$).}}
(maksimalno) rje\v senje na malom intervalu
$(x_0-\varepsilon,x_0+\varepsilon)$ oko $x_0$.
\end{teorem}
{\bf Posebni slu\v cajevi.}
Dva posebna slu\v caja Picard-Lindel\"of-ove teoreme \' ce biti od posebnog
zna\v caja za nas:
\begin{enumerate}
\item ako je $y\mapsto f(y)$ diferencijabilna onda $(\star)$
   ima jedinstveno lokalno rje\v senje;
\item ako je $y\mapsto f(x,y)=A(x)y$ linearna onda $(\star)$
   ima jedinstveno globalno (!) rje\v senje $y:I\to\R^n$.
\end{enumerate}
\begin{dokaz}
Vje\v zba \ref{vjezba:Picard}.
\end{dokaz}

{\bf Potsjetnik iz analiti\v cke geometrije}
{\it Ravan} u koordinatnom prostoru mo\v ze biti definisana na vi\v se na\v cina. Ravan je data u {\it parametarskoj formi} ukoliko je pretstavljena kao ta\v cka u prostoru koja ima neku linearnu funkciju kao parametar ili dvije parametarske varijable (recimo $t$ i $s$) u svakoj od svojih koordinata. Na primjer, $(1+t-5s,s+t,5-7t), t, s \in \mathbb{R},$ je parametarski data ravan u $\mathbb{R}^3$.

Koriste\' ci sli\v can pristup kao kod parametarski datih pravih u $\mathbb{R}^3$, ukoliko su nam date $3$ ta\v cke na ravni
$A$, $B$ i $C=$, {\it koje nisu kolinearne}, prvo defini\v semo vektore $$u=(x_b-x_a,y_b-y_a,z_b-z_a)\ i \ v=(x_c-x_a,y_c-y_a,z_c-z_a)$$ koji idu od A do B i od A do C.
Sada \' ce svi vektori oblika
$(x_a,y_a,z_a)+t\cdot u+s\cdot v$ za sve vrijednosti parametarskih varijabli $t$ i $s$ i\' ci kroz ta\v cke na ravni koja prolazi kroz ta\v cke $A,B,C$.
Stoga ova ravan ima parametarsku formu:
\begin{figure}[t]
\centering
\includegraphics[height=4cm, bb= 0 0 678 550]{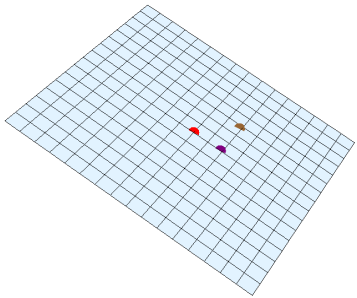}
\caption{Ravan $x+z+\sqrt{2}=0$, sa pripadnim ta\v ckama}\label{ravanSlika}
\end{figure}
\begin{eqnarray*}
x & = & x_a+t(x_b-x_a)+s(x_c-x_a), \\
y & = & y_a+t(y_b-y_a)+s(y_c-y_a), \\
z &=& z_a+t(z_b-z_a)+s(z_c-z_a).
\end{eqnarray*}
\begin{primjer}
Parametarski predstaviti ravan datu na slici \ref{ravanSlika}:
$$
x+z-\sqrt{2} = 0.
$$
Uzmimo tri nekolinearne ta\v cke na ravni $$A\left(\frac{\sqrt{2}}{2},0,\frac{\sqrt{2}}{2}\right), B\left(\frac{\sqrt{2}}{2},1,\frac{\sqrt{2}}{2}\right), C(0,0,\sqrt{2}),$$ kako bismo dobili parametrizaciju ravni
\begin{eqnarray*}
x & = & x_a+\tau(x_b-x_a)+\sigma(x_c-x_a) =  \frac{\sqrt{2}}{2}-\frac{\sqrt{2}}{2}\sigma, \\
y & = & y_a+\tau(y_b-y_a)+\sigma(y_c-y_a) = \tau, \\
z &=& z_a+\tau(z_b-z_a)+\sigma(z_c-z_a) = \frac{\sqrt{2}}{2}+\frac{\sqrt{2}}{2}\sigma.
\end{eqnarray*}
\end{primjer}
\chapter{Du\v zina luka - predznanje}\label{appendixDuzinaLuka}
Ovaj dodatak prati ekspoziciju iz \cite{dedagic2005matematicka}.

\emph{U punoj verziji teksta.}

\chapter{Krivolinijski integral prve vrste - predznanje}
\vspace{1cm}
Ovaj dio prati ekspoziciju iz \cite{okicic2014matematika}.

\emph{U punoj verziji teksta.}


%

\pagestyle{uvod}
\bibliography{difgeo}{}
\bibliographystyle{plain}
\pagestyle{uvod}

\end{document}